\documentclass[a4paper,12pt]{amsart}

\usepackage{url}
\usepackage[margin=0.9in]{geometry}
\usepackage[english]{babel}
\usepackage[utf8]{inputenc}

\usepackage{amsmath}
\usepackage{amssymb}
\usepackage{mathbbol}
\usepackage{algorithm}
\usepackage{algorithmic}
\usepackage{cite}
\usepackage{graphicx}
\usepackage{xcolor}
\usepackage{tikz}

\newtheorem{theorem}{Theorem}[section]

\theoremstyle{definition}
\newtheorem{definition}[theorem]{Definition}
\newtheorem{example}[theorem]{Example}

\theoremstyle{remark}
\newtheorem{remark}[theorem]{Remark}

\numberwithin{equation}{section}

\begin{document}

\title[On the computation of modular forms on noncongruence subgroups]{On the computation of modular forms on noncongruence subgroups}


\author[D. Berghaus]{David Berghaus}
\address{Bethe Center, University of Bonn, Nussallee 12, 53844 Bonn, Germany}
\curraddr{}
\email{berghaus@th.physik.uni-bonn.de}
\thanks{}

\author[H. Monien]{Hartmut Monien}
\address{Bethe Center, University of Bonn, Nussallee 12, 53844 Bonn, Germany}
\curraddr{}
\email{hmonien@uni-bonn.de}
\thanks{}

\author[D. Radchenko]{Danylo Radchenko}
\address{Laboratoire Paul Painlevé, University of Lille, F-59655 Villeneuve d'Ascq, France}
\curraddr{}
\email{danradchenko@gmail.com}
\thanks{}


\date{}

\dedicatory{}

\newcommand{\Hbar}{\overline{\mathcal{H}}}
\newcommand{\bigO}{\mathcal{O}}
\newcommand{\ZZ}{\mathbb{Z}}
\newcommand{\QQ}{\mathbb{Q}}
\newcommand{\RR}{\mathbb{R}}
\newcommand{\CC}{\mathbb{C}}
\newcommand{\subgroup}{\leqslant}
\newcommand{\epsmachine}{\epsilon_\textrm{machine}}
\newcommand{\iu}{{i\mkern1mu}}

\begin{abstract}
	We present two approaches that can be used to compute modular forms on 
	noncongruence subgroups. The first approach uses Hejhal's method for which we 
	improve the arbitrary precision solving techniques so that the algorithm becomes 
	about up to two orders of magnitude faster in practical 
	computations. This allows us to obtain high precision numerical estimates of the 
	Fourier coefficients from which the algebraic expressions can be identified
	using the LLL algorithm. The second approach is restricted to genus zero 
	subgroups and uses efficient methods to compute the Belyi map from which the 
	modular forms can be constructed.
\end{abstract}

\maketitle
\allowdisplaybreaks

\section{Introduction}

Congruence subgroups of the modular group play a significant role in
number theory and have been studied extensively. On the the other hand
noncongruence subgroups and their modular forms are still poorly
understood although some progress has been achieved recently by Chen
\cite{MR3788845} providing a moduli interpretation of noncongruence modular
curves and Calegari, Dimitrov and Tang proving the unbounded
denominator conjecture \cite{https://doi.org/10.48550/arxiv.2109.09040}.

Still the efficient computation of modular forms on noncongruence subgroups of the
modular group remains an open problem due to the lack of non-trivial Hecke operators \cite{berger,li_long_yang}.
The computations of the coefficients of the Fourier expansions of
noncongruence modular forms have therefore so far typically been
limited to special types of subgroups such as noncongruence character
groups \cite{kurt_long} and examples of low number field degree and
index \cite{asd, fiori_franc,MR2763946}. 
Recent advances have been made by the second author who computed the Hauptmodul for a few genus zero subgroups of large index \cite{monien_j2,monien_co3}.

The aim of this paper is to present effective numerical methods in
order to obtain more data on modular forms of noncongruence subgroups
in a systematic way. The outline of the paper is as follows: Section \ref{sec:background} provides the necessary mathematical background and notation, Section \ref{sec:hejhal} describes a numerical method to compute Fourier coefficients of modular forms for general subgroups that is due to Hejhal \cite{hejhal1999} and uses modular transformations to obtain a linear system of equations that can be solved to obtain approximations of the Fourier coefficients of modular forms of arbitrary weight. While Hejhal's method is very versatile, its limitation in practical computations has been that the linear solving involved becomes very slow when applied to high precision. To overcome this difficulty, we demonstrate in Section \ref{sec:noncong_numerics} that mixed-precision iterative solving techniques can be used to significantly improve the performance of Hejhal's method, making the computation of examples that have previously been out of reach feasible. Finally, in Section \ref{sec:genus_zero_belyi_maps_and_forms} we present an alternative approach that is restricted to genus zero subgroups. For this approach we make use of efficient methods to compute genus zero Belyi maps and demonstrate how Fourier expansions of modular forms can be obtained from these.

\section{Background and notation}
\label{sec:background}

Let $\mathrm{SL}(2,\ZZ)$ denote the group of all integer $2 \times 2$ matrices with
determinant 1.  An element
\begin{equation}
  \gamma = \begin{pmatrix}
    a & b\\
    c & d
  \end{pmatrix} \in \mathrm{SL}(2,\ZZ)\,,
\end{equation}
acts on the upper half plane
$\mathcal{H} := \{\tau \in\mathbb{C}\, | \, \textrm{Im}(\tau)>0\}$
in a standard way via Möbius transformations
\begin{equation}
	\gamma(\tau) := \frac{a\tau+b}{c\tau+d}\,.
\end{equation}
Note that
\begin{equation}
	\textrm{Im}(\gamma(\tau)) = \frac{\textrm{Im}(\tau)}{|c\tau+d|^2} > 0\,,
\end{equation}
which means that the elements $\gamma(\tau)$ are also on the upper half plane. It is also immediate to see that $\gamma$ and $-\gamma$ act in the same way. For this
reason, it is often more natural to work with the projective group
\begin{equation}
	\mathrm{PSL}(2,\ZZ) \simeq \mathrm{SL}(2,\ZZ)/\{\pm \mathbb{1}\}\,.
\end{equation}
In the following, we denote $\mathrm{PSL}(2,\ZZ)$ by $\Gamma$ and refer to it as the
modular group.
\begin{definition}[Modular Form]
	Let $f(\tau)$ be a holomorphic function from $\mathcal{H}$ to $\mathbb{C}$. Let $G \subgroup \Gamma$ be a finite index subgroup of $\Gamma$. Then we say that $f(\tau)$ is a modular form on $G$ if it satisfies the functional equation
	\begin{equation}
		f(\gamma(\tau)) = (c\tau+d)^k f(\tau)\,,
	\end{equation}
	for all $\gamma$ in $G$.
\end{definition}
The number $k \in 2\mathbb{N}$ is called the weight of $f$ and $(c\tau+d)^k$ is
the so-called automorphy factor. (More general definitions of modular forms including
odd weights and multiplier system exist but we will not consider them in this work.)
Furthermore, we say that a modular form $f$ is \cite{cohen_stroemberg2017}:
\begin{enumerate}
	\item \emph{weakly holomorphic} if $f$ is holomorphic in $\mathcal{H}$ but might have poles at the boundary $\partial \mathcal{H} := \mathbb{Q} \cup \{\iu\infty\}$.
	\item \emph{holomorphic} if $f$ is holomorphic in $\Hbar := \mathcal{H} \cup \partial \mathcal{H}$.
	\item a \emph{cusp form} if $f$ vanishes at $\partial \mathcal{H}$.
\end{enumerate}
In what follows, when we say modular form we will usually mean holomorphic modular form. Additionally, weakly holomorphic modular forms of weight zero are often called
\emph{modular functions}. It is convenient to introduce the \emph{slash operator}
\begin{equation}
	(f|_k\gamma)(\tau) := (c\tau +d)^{-k}f(\gamma(\tau))\,,
\end{equation}
which defines a right action of $\Gamma$ on the space of complex-valued functions,
i.e.,
\begin{equation}
	f|_k \gamma_1 |_k \gamma_2 = f|_k \gamma_1 \gamma_2\,.
\end{equation}

\subsection{Fundamental domains}
\label{sec:fund_domains}
We define a fundamental domain of a group $G \subgroup \Gamma$ as follows:
\begin{definition}[Fundamental Domain {\cite[Definition 4.3.1]{cohen_stroemberg2017}}]
	A closed set $\mathcal{F}(G) \subset \Hbar$ is said to be a \emph{fundamental domain} if
	\begin{enumerate}
		\item For any point $\tau \in \Hbar$ there is a $\gamma \in G$ such that $\gamma(\tau) \in \mathcal{F}(G)$.
		\item If for any points $\tau$ and $\tau' := \gamma(\tau)$ we have $\tau \neq \tau'$ then $\tau, \tau' \in \partial\mathcal{F}(G)$.
	\end{enumerate}
\end{definition}
Note that $\Gamma$ can be generated by the elements
\begin{equation}
	\label{eq:sl2z_generators}
	S = \begin{pmatrix}
		0 & -1\\
		1 & 0
	\end{pmatrix}
	\quad
	\textrm{and}
	\quad
	T = \begin{pmatrix}
		1 & 1\\
		0 & 1
	\end{pmatrix}\,.
\end{equation}
The matrix $S$ therefore corresponds to the action $\tau \rightarrow -1/\tau$ which can be
viewed as an inversion while $T$ corresponds to the action $\tau
\rightarrow \tau+1$, a translation. Moreover, we have the relations
\begin{equation}
	S^2 = \mathbb{1} \quad \textrm{and} \quad (ST)^3 = \mathbb{1}\,.
\end{equation}
A fundamental domain for the modular group is given by the set
\begin{equation}
	\label{eq:modular_group_fund_domain}
	\mathcal{F}(\Gamma) = \{\tau \in \Hbar \,,\, |\tau| \geq 1 \, \textrm{and} \, |\textrm{Re}(\tau)| \leq 1/2\} \cup \{\iu\infty\}\,.
\end{equation}
The fundamental domain $\mathcal{F}(\Gamma)$ has three points that play a special role:
\begin{enumerate}
	\item $\iu\infty$: a \emph{cusp};
	\item $\iu$: an \emph{elliptic point of order 2} which has a non-trivial
	stabilizer $S$ with $S^2 = \mathbb{1}$;
	\item $\rho = \exp(2\pi \iu/3)$: an \emph{elliptic point of order 3} which has a
	non-trivial stabilizer $ST$ with $(ST)^3 = \mathbb{1}$ (alternatively we could
	also choose the point $-\bar{\rho} = \exp(\pi \iu/3)$).
\end{enumerate}
Moreover, since $\gamma(\iu\infty) = a/c$, we can see that the cusps are located at
$\mathbb{P}^1(\QQ) = \{\iu\infty\} \cup \QQ$.
For a finite index subgroup $G \subgroup \Gamma$ of index $\mu$, a fundamental domain for $G\setminus \mathcal{H}$ is given by
\begin{equation}
	\mathcal{F}(G) = \cup_{i=1}^\mu \gamma_i \mathcal{F}(\Gamma)\,,
\end{equation}
where $\gamma_i$ are right coset representatives of $G\setminus\Gamma$. The suitably defined quotient
$G\setminus \Hbar$ (see, e.g., \cite[Theorem 4.4.3]{cohen_stroemberg2017}) is a Riemann surface whose genus can be
computed using the
formula \cite[Proposition 5.6.17]{cohen_stroemberg2017}
\begin{equation}
	g = 1 + \frac{\mu}{12} - \frac{n(e_2)}{4} - \frac{n(e_3)}{3} - \frac{n(c)}{2}\,,
\end{equation}
where $n(e_2)$, $n(e_3)$ denote the amount of inequivalent elliptic points of order two and three, respectively, and $n(c)$ denotes the amount of cusp representatives.
\begin{definition}[Signature]
	\label{def:signature}
	We define the \emph{signature} of $G\subgroup \Gamma$ to be the tuple $(\mu,g,n(c),n(e_2),n(e_3))$. Note that a signature does not uniquely specify $G$!
\end{definition}
We call the maps $A_j \in \textrm{PSL}(2,\ZZ)$ that map $\iu\infty$ to the cusp $p_j$
on the real line
\begin{equation}
	\label{eq:cusp_normalizer}
	A_j(\iu\infty) = p_j\,,
\end{equation}
and satisfy
\begin{equation}
	A_j^{-1} S_j = T^N\,,
\end{equation}
the \emph{cusp normalizers}, where $S_j$ is the generator of the stabilizer of $p_j$
(we use the notation of Strömberg \cite{stroemberg_recent}, some authors use the
reversed notation) and $N$ denotes the cusp width at infinity.

\subsection{Subgroups of the modular group}
Let $N$ be a positive integer. Then we call
\begin{equation}
	\Gamma(N) := \left\{\begin{pmatrix}
		a & b\\
		c & d
	\end{pmatrix}
	\equiv
	\begin{pmatrix}
		1 & 0\\
		0 & 1
	\end{pmatrix}
	\pmod{N}
	\quad \textrm{and} \quad
	\begin{pmatrix}
		a & b\\
		c & d
	\end{pmatrix}
	\in \Gamma
	\right\}\,,
\end{equation}
the \emph{principal congruence subgroup of level N}. The index of $\Gamma(N)$ is given by \cite[Corollary 6.2.13]{cohen_stroemberg2017}
\begin{equation}
	\left[\Gamma:\Gamma(N)\right] = \frac{1}{2} N^3
	\prod_{p|N}\left(1-\frac{1}{p^2}\right)\,.
\end{equation}
\begin{definition}[Congruence Subgroup]
	A subgroup $G \subgroup \Gamma$ is a \emph{congruence subgroup} of level $N$ iff
	it contains $\Gamma(N)$ for some $N\in \ZZ^+$ (i.e., if $\Gamma(N) \subgroup G$
	for some $N$).
\end{definition}
Important examples of congruence subgroups are
\begin{equation}
	\label{eq:Gamma_0_N}
	\Gamma_0(N) := \left\{\begin{pmatrix}
		a & b\\
		c & d
	\end{pmatrix}
	\equiv
	\begin{pmatrix}
		* & *\\
		0 & *
	\end{pmatrix}
	\pmod{N}
	\quad \textrm{and} \quad
	\begin{pmatrix}
		a & b\\
		c & d
	\end{pmatrix}
	\in \Gamma
	\right\}\,,
\end{equation}
and
\begin{equation}
	\Gamma_1(N) := \left\{\begin{pmatrix}
		a & b\\
		c & d
	\end{pmatrix}
	\equiv
	\begin{pmatrix}
		1 & *\\
		0 & 1
	\end{pmatrix}
	\pmod{N}
	\quad \textrm{and} \quad
	\begin{pmatrix}
		a & b\\
		c & d
	\end{pmatrix}
	\in \Gamma
	\right\}\,,
\end{equation}
which satisfy
\begin{equation}
	\Gamma(N) \subgroup \Gamma_1(N) \subgroup \Gamma_0(N) \subgroup \Gamma\,.
\end{equation}
Subgroups that are not congruence are called \emph{noncongruence subgroups}. It has
been proven by Stothers \cite{stothers_1984} that noncongruence subgroups are much
more numerous than congruence subgroups (in the sense that the proportion of the
latter among all subgroups of index $n$ goes to $0$ as $n\to\infty$). An algorithm to
test if a given group $G$ is congruence or not has been given by Hsu \cite{hsu}.

A useful tool for studying subgroups $G \subgroup \Gamma$ is the interpretation of
the action of $G$ on the cosets of $G\setminus \Gamma$ as an action of the symmetric
group $S_\mu$. This theory has been developed by Millington \cite{millington} and its
usefulness when performing computations with subgroups of $\Gamma$ has first been
demonstrated by Atkin and Swinnerton-Dyer \cite{asd}.
\begin{definition}[Legitimate Pair, \cite{millington}]
	A pair $(\sigma_S,\sigma_R)$ with $\sigma_S,\sigma_R \in S_\mu$ is called
	\emph{legitimate} if $\sigma_S^2=\sigma_R^3=\mathbb{1}$ and if the group $\Sigma$
	that is generated by $\sigma_S$ and $\sigma_R$ is transitive.
\end{definition}
\begin{definition}[Equivalence Modulo 1, \cite{millington}]
	Two legitimate pairs $(\sigma_S,\sigma_R)$ and $(\sigma'_S,\sigma'_R)$ are said
	to be \emph{equivalent (modulo 1)} if there exists a
	$\sigma \in S_\mu$ such that
	$(\sigma^{-1}\sigma'_S\sigma,\sigma^{-1}\sigma'_R\sigma) = (\sigma_S,\sigma_R)$
	and $\sigma(1) = 1$ (i.e., that $\sigma$ fixes 1).
\end{definition}
\begin{theorem}[Millington]
	\label{th:millington}
	There is a one-to-one correspondence between subgroups $G$ of index $\mu$ in
	$\Gamma$ and equivalence classes modulo 1 of legitimate pairs
	$(\sigma_S,\sigma_R)$. Moreover, $n(e_2)$ and $n(e_3)$ are given by the number
	of fixed elements of $\sigma_S$ and $\sigma_R$, respectively, and $n(c)$
	is the number of elements that are fixed by $\sigma_T = \sigma_S
	\sigma_R$. Additionally, the cycle structure of $\sigma_T$ reflects the cusp
	widths of $G$.
	\begin{proof}
		See~\cite[Theorem 2]{millington}
	\end{proof}
\end{theorem}
The action of $\Gamma$ on the cosets of $G$ gives rise to a map
\begin{equation}
	\label{eq:millington_homo}
	\phi \,:\, \Gamma \rightarrow S_\mu\,,
\end{equation}
which satifies $\phi(x\cdot y) = \phi(x)\cdot \phi(y)$ and is hence a homomorphism. Note that the set of
coset representatives
$\gamma_i$, $i=1,...,\mu$ of $G$ satisfies $\phi(\gamma_i)(1) = i$ (see
\cite{stroemberg_recent} for more details).

Millington's theorem also provides a method to list all subgroups of a given index by
filtering legitimate pairs into equivalence classes modulo 1. This algorithm has been
applied by Strömberg \cite{stroemberg_recent} to calculate representatives of all
subgroups in $\Gamma$ with $\mu \leq 17$ up to relabelling (or in other words,
conjugation in $\Gamma$). Strömberg has released this data in \cite{stroemberg_db}.
\begin{example}[$\Gamma_0(5)$]
  \label{ex:Gamma_0_5}
  Consider the group $\Gamma_0(5)$ (defined as in
  Eq. \eqref{eq:Gamma_0_N}) with signature $(6, 0, 2, 2, 0)$. As
  a legitimate pair for $\Gamma_0(5)$ one can choose $\sigma_S =
  (1)(2)(3\,4)(5\,6)$ and $\sigma_R =
  (1\,2\,3)(4\,5\,6)$. Following from this, we get that
  $\sigma_T = \sigma_S \sigma_R = (1\,2\,3\,5\,4)(6)$.
  A set of right coset representatives can be chosen to be
  \begin{equation}
    \left\{
      \begin{pmatrix}
        1 & 0\\
        0 & 1
      \end{pmatrix},\,
      \begin{pmatrix}
        1 & 1\\
        0 & 1
      \end{pmatrix},\,
      \begin{pmatrix}
        1 & 2\\
        0 & 1
      \end{pmatrix},\,
      \begin{pmatrix}
        1 & -1\\
        0 & 1
      \end{pmatrix},\,
      \begin{pmatrix}
        1 & -2\\
        0 & 1
      \end{pmatrix},\,
      \begin{pmatrix}
        -2 & -1\\
        1 & 0
      \end{pmatrix}
    \right\}\,,
  \end{equation}
  which can be expressed as words in $S$ and $T$ as follows
  \begin{equation}
    \left\{
      \mathbb{1},\,
      T,\,
      T^2,\,
      T^{-1},\,
      T^{-2},\,
      T^{-2}S
    \right\}\,.
  \end{equation}
  A fundamental domain and the corresponding coset labels can therefore be chosen
  as in Fig.~\ref{fig:Gamma0_5_fund_domain_labeled}. We can see that this group has
  two cusps: One of width 5 at $\iu\infty$ and one of width 1 at $-2$.
  Additionally, we can tell from the signature and by looking at $\sigma_R$ that
  $\Gamma_0(5)$ has no elliptic points of order three. The two elliptic points of
  order two are located at $\gamma_1(\iu)$ and $\gamma_2(\iu)$ where $\gamma_j$
  corresponds to the coset representative of label $j$ because 1 and 2 are fixed by
  $\sigma_S$.
\end{example}
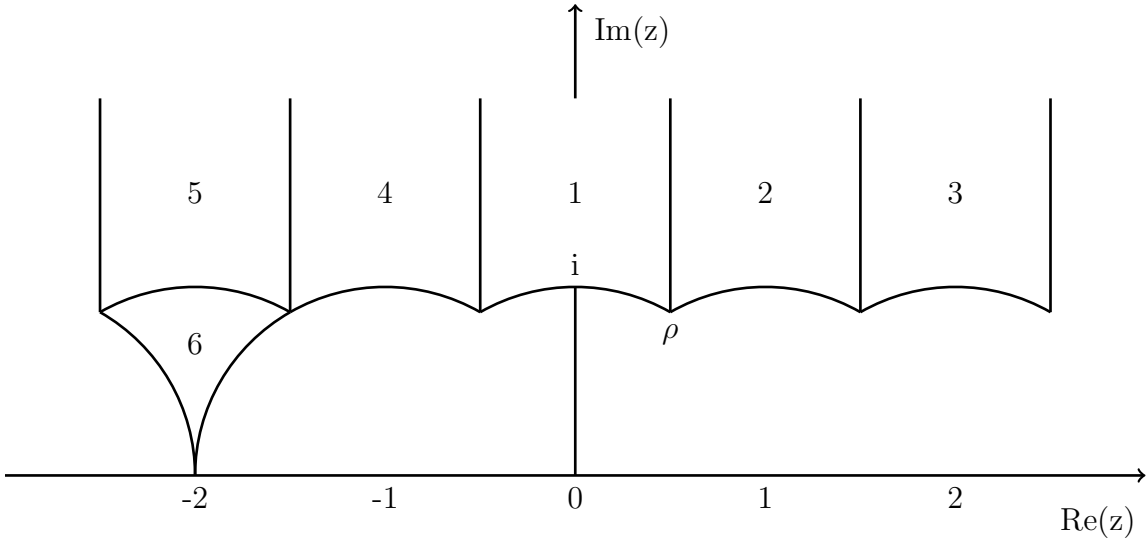
\begin{figure}[h]
  \centering
  \begin{tikzpicture}[line width=1pt,scale=2.5]
    \draw[->] (-3,0) -- (3,0); 
    \foreach \x in {-2,...,2} \node[below] at (\x,0) {\x};
    \node[below] at (2.75,-0.1) {$\rm{Re}(z)$};
    \draw (0,0) -- (0,1);
    \draw[->] (0,2) -- (0,2.5);
    \node[right,below] at (0.3,2.5) {$\rm{Im}(z)$};
    \node[above] at (0,1) {i};
    \node[below] at ({1/2},{sqrt(3)/2)}) {$\rho$};
   \foreach \center / \name in {(-2,1.5)/5, (-1,1.5)/4, (0,1.5)/1,
     (1,1.5)/2, (2,1.5)/3, (-2, 0.7)/6} \node at \center {\name};
    \foreach \x in {-2,...,2} \draw ({\x-1/2},{sqrt(3)/2)}) arc (90+30:90-30:1);
    \foreach \x in {-2,...,3} \draw ({\x-1/2},{sqrt(3)/2)}) -- ({\x-1/2},2);
    \draw (-2,0) arc (0:60:1);
    \draw (-2,0) arc (180:120:1);
  \end{tikzpicture}
  \caption{A fundamental domain for $\Gamma_0(5)$ corresponding to the legitimate pair  $\sigma_S = (1)(2)(3\,4)(5\,6)$ and $\sigma_R = (1\,2\,3)(4\,5\,6)$.}
  \label{fig:Gamma0_5_fund_domain_labeled}
\end{figure}

\subsection{Fourier expansions of modular forms}
We have seen in the previous sections that modular forms are functions on the upper
half plane that satisfy certain functional equations. Additionally we have seen that
the cusp widths are always finite and that the modular forms are periodic with
respect to these cusp widths. Modular forms can therefore be expanded
as Fourier series in variable
\begin{equation}
	q_N := \exp(2\pi \iu \tau/N) = \exp(2\pi \iu (x+iy)/N) = \exp(2\pi \iu x/N)\exp(-2\pi y/N)\,,
\end{equation}
where $N$ denotes the cusp width and $\tau = x+\iu y \in \mathcal{H}$ (we will also
often use the convention $q := q_1$). It is important to note that $q_N$ decays
exponentially as $y\rightarrow \infty$. If $f$ is a modular form and the cusp
width at $\iu\infty$ is given by $N$, then we can write
\begin{equation}
	f(\tau) = \sum_{n=-\infty}^{\infty} a_n q_N^n\,,
\end{equation}
with $a_i \in \CC$. For congruence subgroups, it is known that there exist bases of modular forms whose Fourier coefficients are defined over $\QQ$ or cyclotomic fields. For noncongruence subgroups the
Fourier coefficients
are defined over $\bar{\QQ}$ and are of the form (see for example
Atkin-Swinnerton-Dyer \cite{asd})
\begin{equation}
	a_n = u^m b_n\,,
\end{equation}
where $b_n$ and $u^N$ are defined over a number field $K$ which is generated over $\QQ$ by an algebraic number $v$ (i.e., $K = \QQ(v)$).
\begin{definition}[Valuation of a modular form]
	We define the valuation of a modular form to be the index of the first non-zero
	Fourier coefficient.
\end{definition}
\begin{remark}
	By using the valuation of a modular form, many properties immediately follow from its $q_N$ expansion. For example, a modular form can only be holomorphic if its Fourier expansion starts at $n \geq 0$ because negative values of $n$ would lead to poles at $\iu\infty$ due to the decay of $q_N$. Following the same argument, cusp forms need to have Fourier expansions that start with $n>0$.
\end{remark}

\subsection{Spaces of modular forms}
We denote the space of holomorphic modular forms of (even) weight $k$ on $G$ by $M_k(G)$ and similarly define $S_k(G)$ to be the space of cusp forms. The dimensions of these spaces can be computed from their signature \cite[Theorem 5.6.18]{cohen_stroemberg2017}
\begin{align}
	\textrm{dim}(M_k(G)) &= (k-1)(g-1) + \left\lfloor\frac{k}{4}\right\rfloor n(e_2) + \left\lfloor\frac{k}{3}\right\rfloor n(e_3) + \left\lfloor\frac{k}{2}\right\rfloor n(c)\,,\\
	\textrm{dim}(S_k(G)) &= \textrm{dim}(M_k(G)) - n(c) + \delta_{k,2}\,.
\end{align}

\subsection{Hauptmoduls}
\label{sec:hauptmoduls}
Subgroups $G \subgroup \Gamma$ of genus zero have a special type of modular function
called the \emph{Hauptmodul} (which we denote by $j_G$).
\begin{definition}[Hauptmodul]
	Let $G$ be a subgroup of genus zero. Then a Hauptmodul is any isomorphism
	\begin{equation}
		j_G\, : \, G\backslash\mathcal{\Hbar}\,\rightarrow\, P^1(\CC)\,.
	\end{equation}
\end{definition}
Because the modular group $\Gamma$ has genus zero, it has a Hauptmodul which is
referred to as \emph{Klein invariant} or \emph{modular j-invariant}. Its Fourier
expansion is given by
\begin{equation}
	j(\tau) = \frac{E_4^3}{\Delta(\tau)} = q^{-1} + 744 + 196884q + 21493760q^2 + 864299970q^3 + ...\,.
\end{equation}
and its values at the elliptic points are
\begin{equation}
	j(\iu) = 1728 \quad \textrm{and} \quad j(\rho) = 0\,.
\end{equation}
Because $j$ has negative valuation, one can also see that is has a pole of order 1 at infinity.
We remark that the Hauptmodul can be chosen uniquely up to a constant term. The choice of $744$ for the constant term has historical reasons. For groups $G \neq \Gamma$ we will instead set the constant term to zero and use the normalization
\begin{equation}
	j_G(\tau) = q_N^{-1} + 0 + \sum_{n=1} a_n q_N^n\,,
\end{equation}
which specifies $j_G$ uniquely \cite{asd}.
\begin{theorem}
	Let $f$ be a meromorphic function on $\mathcal{H}$. The following statements are equivalent:
	\begin{enumerate}
		\item $f$ is a modular function for $\Gamma$ of weight 0.
		\item $f$ is a quotient of two modular forms for $\Gamma$ of equal weight.
		\item $f$ is a rational function of $j$.
	\end{enumerate}
	\begin{proof}
		See \cite[Theorem 5.7.3]{cohen_stroemberg2017}
	\end{proof}
\end{theorem}
\begin{theorem}
	\label{theorem:holomorphic_modform_polynomial}
	Every modular function on $G$ that is holomorphic outside $\iu\infty$ can be written as a polynomial $P(j_G(\tau))$.
	\begin{proof}
		See Cox \cite[Lemma 11.10 (ii)]{cox} for the case of $G = \Gamma$ (the proof for general $G$ is analogous).
	\end{proof}
\end{theorem}

\section{Hejhal's method}
\label{sec:hejhal}
A general method to compute numerical approximations of the coefficients of modular
forms in an expansion basis has been given by Hejhal \cite{hejhal1999} (based on an
idea of Stark) who has developed this method to compute Maass cusp forms on Hecke
triangle groups. The basic idea of Hejhal's method is to expand a modular form (for
example in a $q$-expansion basis) and to afterwards impose the modular transformation
property of the expansion on a finite set of points. This creates a linear system of
equations that can be solved to obtain numerical approximations of the expansion
coefficients. Due to the generality of this method (in principle the only
requirements are a converging expansion basis for the modular form and a automorphy
condition) it has since then been adapted by many authors. For example, Selander and
Strömbergsson \cite{selander-stroembergsson} generalized the method for fundamental
domains with multiple cusps to compute some examples of genus 2 coverings and
Strömberg used this method to compute Maass cusp forms for $\Gamma_0(N)$ and
non-trivial multiplier systems \cite{stroemberg_gamma_0_N} as well as Maass cusp forms
for noncongruence subgroups \cite{stroemberg_recent}. Applications of Hejhal's method
using arbitrary precision arithmetic have been performed by Booker, Strömbergsson and
Venkatesh \cite{mass_multiprecision} who computed the first ten Maass cusp forms of
$\Gamma$ to 1000 digits precision, Bruinier and Strömberg \cite{harmonic_weak} who
computed harmonic weak Maass cusp forms and Voight and Willis \cite{voight_willis}
(see also the improved method in KMSV \cite{klug_musty_schiavone_voight_2014}) who
computed Taylor expansions of modular forms.

\subsection{The case $G=\Gamma$}
\label{sec:hejhal_single_cusp}
To illustrate Hejhal's method \cite{hejhal1999} we first consider the simplest case
$G=\Gamma$ for which the fundamental domain only has a single cusp and whose
fundamental domain is given by Eq. \eqref{eq:modular_group_fund_domain}. The point
inside $\mathcal{F}(\Gamma)$ with the smallest height (i.e., the smallest imaginary
value) is given by $\rho$ whose height is $Y_0 = \sqrt{3}/2$. Now we choose a set of
$2Q$ points $\tau_m$ that are equally spaced between $-1/2$ and $1/2$ along a
horizontal line with height $Y < Y_0$
\begin{equation}
	\label{eq:horocycle_points}
	\tau_m = x_m + iY = \frac{1}{2Q}\left(m-Q+\frac{1}{2}\right) + iY, \quad 0 \leq m \leq 2Q-1,\quad Y < Y_0\,.
\end{equation}
\begin{remark}
	We will always choose $Y = 0.8\cdot Y_0$ throughout this paper.
\end{remark}
We also refer to the points $\tau_m$ on this horizontal line as a
\emph{horocycle}. Note that because these
points are located \emph{below} $\mathcal{F}(\Gamma)$, they are all \emph{outside}
$\mathcal{F}(\Gamma)$. Now for each point $\tau_m$ there exists a map $\gamma_m \in
\Gamma$ such that
\begin{equation}
	\tau_m^* = \gamma_m(\tau_m) \in \mathcal{F}(\Gamma), \quad \gamma_m = \begin{pmatrix}
		a_m & b_m\\
		c_m & d_m
	\end{pmatrix} \in \Gamma.
\end{equation}
We call the maps $\gamma_m$ the \emph{pullback} to the fundamental domain. In the
case of the modular group, finding such a pullback map is straightforward, we simply
need to form words in the generators $S \rightarrow -1/\tau$ and $T \rightarrow
\tau+1$ depending on if $|\tau| < 1$, $\textrm{Re}(\tau) < -1/2$ or
$\textrm{Re}(\tau) > 1/2$ and form the matrix products. Afterwards, we expand the
modular form in a suitable basis (which is, in our case, given by powers of $q$) up
to a finite order $M_0 := M(Y_0)$ so that our expansion converges inside
$\mathcal{F}(\Gamma)$ up to the machine epsilon $\epsmachine$. The value of $M_0$ can
be guessed in advance by using the asymptotic growth conditions of the coefficients
(the coefficients of cusp forms have asymptotic growth $\bigO(n^{k/2})$ and for
holomorphic modular forms the coefficients grow like $\bigO(n^k)$ \cite{serre}).
Although such a choice of $M_0$ works well in practice, it is non-rigorous and there
is therefore no guarantee at this point that the result will be correct. This is one
of the reasons why it is difficult to make Hejhal's method rigorous. In order to be a
modular form, the expansion now needs to (at least numerically) match the automorphy
condition
\begin{equation}
	\label{eq:hejhal_original}
	f(\tau_m) \approx \sum_{n=0}^{M_0}a_n {q(\tau_m)}^n \overset{!}{=} (c_m\cdot \tau_m+d_m)^{-k} f(\tau_m^*) \approx (c_m\cdot \tau_m+d_m)^{-k} \sum_{n=0}^{M_0}a_n {q(\tau_m^*)}^n\,,
\end{equation}
where $q(\tau) = \exp(2\pi i \tau)$ (we illustrate this method here for the example of holomorphic modular forms, but it can obviously be applied analogously for cusp forms or Hauptmoduls). For numerical reasons it is preferable to work with
\begin{equation}
	F(\tau) = y^{k/2}f(\tau)\,,
\end{equation}
where $y = \textrm{Im}(\tau)$.
The function $F$ transforms like
\begin{equation}
	F(\tau_m) = \frac{|c_m\cdot \tau_m+d_m|^k}{(c_m\cdot \tau_m+d_m)^k}F(\tau_m^*)\,,
\end{equation}
and its automorphy factor hence does not change the order of magnitude. Eq. \eqref{eq:hejhal_original} creates a linear system of equations that can in principle be solved to obtain numerical approximations of the expansion coefficients (see for example \cite{asd,hejhal1992}). From a numerical analysis perspective, the resulting linear system of equations however typically becomes ill-conditioned. A more numerically stable approach has been given by Hejhal in \cite{hejhal1999} and uses the Fourier integral formula
\begin{equation}
	a_n Y^{\frac{k}{2}} \exp(-2\pi n Y) = \int_{-\frac{1}{2}}^\frac{1}{2} F(\tau)\exp(-2\pi \iu n x) dx\,,
\end{equation}
where $Y$ denotes the height of the horocycle. Discretizing this integral to approximate it numerically gives
\begin{equation}
	a_n Y^{\frac{k}{2}} \exp(-2\pi n Y) \approx \frac{1}{2Q} \sum_{m=0}^{2Q-1} F(\tau_m) \exp(-2\pi \iu n x_m)\,,
\end{equation}
where $Q > M(Y)$ and $\tau_m$ are again given by Eq.\ \eqref{eq:horocycle_points}. Hejhal then incorporates the automorphy condition by replacing $F(\tau_m)$ with the corresponding pullback
\begin{align}
	a_n Y^{\frac{k}{2}} \exp(-2\pi n Y)
	&\approx \frac{1}{2Q} \sum_{m=0}^{2Q-1} \left( \frac{|c_m \tau_m + d_m|}{\left(c_m \tau_m + d_m\right)}\right)^k F(\tau^*_m)\exp(-2\pi \iu n x_m)\,,\\
	&= \sum_{l=0}^{M_0} a_l \frac{1}{2Q} \sum_{m=0}^{2Q-1} \left( \frac{|c_m \tau_m + d_m|}{\left(c_m \tau_m + d_m\right)}\right)^k (y^*_m)^{\frac{k}{2}} \exp(2\pi \iu (l \tau^*_m - n x_m))\,,\\
	&:= \sum_{l=0}^{M_0} a_l V_{n,l}\,,
\end{align}
where
\begin{equation}
	\label{eq:V_nl}
	V_{n,l} := \frac{1}{2Q} \sum_{m=0}^{2Q-1} \left( \frac{|c_m \tau_m + d_m|}{\left(c_m \tau_m + d_m\right)}\right)^k (y^*_m)^{\frac{k}{2}} \exp(2\pi \iu (l \tau^*_m - n x_m))\,.
\end{equation}
Therefore
\begin{equation}
	0 = \sum_{l=0}^{M_0} a_l \tilde{V}_{n,l}\,,
\end{equation}
with
\begin{equation}
	\label{eq:V_tilde_nl}
	\tilde{V}_{n,l} := V_{n,l} - \delta_{n,l} Y^{\frac{k}{2}} \exp(-2\pi n Y)\,.
\end{equation}
The resulting linear system of equations can be solved numerically for example by imposing a reduced row echelon normalization and dropping the first row of $\tilde{V}_{n,l}$. For example for a one-dimensional space of modular forms we can set $a_0 = 1$ which amounts to solving
\begin{equation}
	\begin{pmatrix}
		\tilde{V}_{1,1} & \dots &  \tilde{V}_{1,M_0}\\
		\vdots & \ddots & \vdots \\
		\tilde{V}_{M_0,1} & \dots & \tilde{V}_{M_0,M_0}  \\
	\end{pmatrix}
	\cdot \left( \begin{array}{c}a_1\\ \vdots \\a_{M_0}\\\end{array} \right)
	= \left( \begin{array}{c}-\tilde{V}_{1,0}\\ \vdots \\-\tilde{V}_{M_0,0}\\\end{array} \right)\,.
\end{equation}
The advantage of this method is that the largest entries of each column are now located on the diagonal. This can be seen in Eq. \eqref{eq:V_tilde_nl}: $V_{n,l}$ depends on the pullbacked points which have a larger imaginary value (and hence smaller $q$-values) than the horocycle points located at height $Y$. For this reason
\begin{equation}
	\label{eq:V_tilde_diagonal_terms_are_largest}
	|Y^{\frac{k}{2}} \exp(-2\pi n Y)| > |V_{n,l}|\,,
\end{equation}
and the largest entries of each column are hence located on the diagonal. This means that the linear system of equations that results from this improved method is significantly better conditioned. The precision of the coefficients depends on the diagonal term in Eq. \eqref{eq:V_tilde_nl}. We can therefore expect the $l$-th coefficient (where $1\leq l \leq M_0$) to be correct to approximately $D-\log_{10}^+\left|\frac{1}{Y^{k/2}\exp(-2\pi l Y)}\right|$ digits precision (this is analogous to Maass cusp forms, see \cite{mass_multiprecision}). The \emph{precision loss} of higher order coefficients can hence be controlled by choosing a smaller value of $Y$ (and following from this a larger value of $Q$).

Once the coefficients $a_l$, $l=0,...,M_0$ have been computed to reasonable accuracy,
approximations of higher coefficients with $l'>M_0$ can be obtained from these
without solving any additional linear systems by using (see \cite{hejhal1999})
\begin{equation}
	\label{eq:hejhal_higher_coeffs}
	a_{l'} = \frac{\sum_{l=0}^{M_0} a_l V_{n,l}}{Y^{\frac{k}{2}}\exp(-2\pi l' Y)}\,,
\end{equation}
where $Y$ is reduced for larger $l'$.
\begin{remark}
	To check the precision of the coefficients computed with Hejhal's method heuristically one can repeat the computation with an independent choice of $Y$. This is especially crucial for Maass cusp forms because it is not a priori clear if the computed solution corresponds to a \emph{true} eigenvalue.
\end{remark}

\subsection{The general case}
The general case, including groups that have multiple cusps, has been worked out by Selander and Strömbergsson \cite{selander-stroembergsson} (see also Strömberg \cite{stroemberg_gamma_0_N,stroemberg_recent}) and follows the same ideas but the resulting expressions are more tedious and the pullback maps more difficult to obtain. If $G$ has multiple cusps then we need to incorporate the Fourier expansions at all cusps in order to obtain convergence in $\mathcal{F}(G)$. Let $j=1,...,n(c)$ label the cusps of $G \subgroup \Gamma$.
\begin{definition}(Width absorbing cusp normalizer)
	\label{def:width_absorbing_normalizer}
	Let $A_j$ denote the cusp normalizer of cusp $j$ as defined in Eq. \eqref{eq:cusp_normalizer}. Let $w_j$ denote the width of cusp $j$. We define the \emph{width absorbing cusp normalizer} of cusp $j$ to be the map $\mathcal{N}_j \in \textrm{PSL}(2,\RR)$ such that
	\begin{equation}
		\mathcal{N}_j(\tau) = A_j(w_j\cdot \tau)\,,
	\end{equation}
	and therefore
	\begin{equation}
		\mathcal{N}_j = A_j \cdot \rho_j = A_j \cdot
		\begin{pmatrix}
			\sqrt{w_j} & 0\\
			0 & 1/\sqrt{w_j}
		\end{pmatrix}\,.
	\end{equation}
\end{definition}
By using width absorbing cusp normalizers, the expansion at the $j$-th cusp is given by
\begin{equation}
	(f|_k\mathcal{N}_j)(\tau) = \sum_{n=0}^{\infty}a^{(j)}q^n\,,
\end{equation}
and can hence always be expanded in $q = q_1$ which is useful and simplifies the expressions.
\begin{definition}[Minimal height of $\mathcal{F}(G)$]
	We define the \emph{minimal height} of $\mathcal{F}(G)$ to be the quantity
	\begin{equation}
		Y_0 := \frac{\sqrt{3}}{2 N_\textrm{max}}\,,
	\end{equation}
	where $N_\textrm{max}$ is the largest cusp width of $G$.
\end{definition}
To compute the pullback of $\tau \notin \mathcal{F}(G)$ into $\mathcal{F}(G)$ we make
use of Millington's theorem (see Theorem \ref{th:millington}). The procedure can be
described as follows:
\begin{enumerate}
	\item Compute the pullback of $\tau$ into $\mathcal{F}(\Gamma)$ which creates a word in $S$, $T$, $T^{-1}$.
	\item Insert the corresponding word in $S$, $T$, $T^{-1}$ into the $\textrm{PSL}(2,\ZZ)$ and $S_\mu$ representations to obtain a
	map $\gamma_\tau \in \Gamma$ and its permutation $\sigma_\tau :=
	\phi(\gamma_\tau) \in S_\mu$.
	\item Let $\sigma_{i} := \phi(\gamma_i) \in S_\mu$ denote the permutation
	representations of the coset representatives. Then the pullback goes into the
	(unique) coset of label $j$ for which $\sigma_\tau(\sigma_j(1)) = 1$.
	\item The pullback into $\mathcal{F}(G)$ is hence given by $\gamma_w = \gamma_j\cdot \gamma_\tau \in \Gamma$.
\end{enumerate}
Once the pullback $w = \gamma_w(\tau)$ into $\mathcal{F}(G)$ has been found, we need
to identify the cusp that is the \emph{closest} to the pullbacked point (in the sense
that its Fourier expansion converges the fastest). This gives rise to a function
(following \cite{selander-stroembergsson,stroemberg_gamma_0_N,stroemberg_recent})
\begin{equation}
	I \,:\, \mathcal{H} \rightarrow \{1,...,n(c)\}\,,
\end{equation}
which returns the cusp label $k$ for which the Fourier expansion at the point $w$
converges the fastest. The complete pullback is therefore given by
\begin{equation}
	\tau^* = \left(\mathcal{N}_{I(w)}^{-1}\cdot \gamma_w\right)(\tau)\,.
\end{equation}
These pullback routines have been contributed by Strömberg to \textsc{Psage} \cite{psage} and have been used in this project as well.

Hejhal's method for multiple cusps can be summarized as follows: For each cusp $j$, we choose a fixed amount of equally spaced points along a horocycle and compute their pullbacks into $\mathcal{F}(G)$. Afterwards we \emph{match} the expansion with the cusp whose Fourier expansion on the pullbacked point converges the fastest. This gives
\begin{equation}
	\tau_{m,j}^* = \left(\mathcal{N}_{I(m,j)}^{-1}\cdot \gamma_w\cdot \mathcal{N}_j \right)(\tau_m) = \begin{pmatrix}
		a_{m,j} & b_{m,j}\\
		c_{m,j} & d_{m,j}
	\end{pmatrix}(\tau_m)\,,
\end{equation}
where $I(m,j) := I(w)$.
In analogy to Section \ref{sec:hejhal_single_cusp} we therefore get
\begin{align}
	a_n^{(j)} Y^\frac{k}{2} \exp(-2\pi n Y) &\approx \frac{1}{2Q} \sum_{m=0}^{2Q-1} (F|_k \mathcal{N}_j)(\tau_m)\exp(-2\pi \iu n x_m)\,,\\
	&= \frac{1}{2Q} \sum_{m=0}^{2Q-1} \left(\frac{|c_{m,j}\tau_m + d_{m,j}|}{(c_{m,j}\tau_m + d_{m,j})}\right)^k (F|_k \mathcal{N}_{I(m,j)})(\tau^*_{m,j})\exp(-2\pi \iu n x_m)\,,\\
	&= \sum_{l=0}^{M_0} a_l^{(I(m,j))} \frac{1}{2Q} \sum_{m=0}^{2Q-1} \left(\frac{|c_{m,j}\tau_m + d_{m,j}|}{(c_{m,j}\tau_m + d_{m,j})}\right)^k (y_{m,j}^*)^\frac{k}{2}\exp(2\pi\iu (l\tau_{m,j}^*-n x_m))\,.
\end{align}
For the analogue of Eq. \eqref{eq:V_nl} we hence get
\begin{equation}
	a^{(j)}_n Y^{\frac{k}{2}} \exp(-2\pi n Y) = \sum_{j'=1}^{\kappa}\sum_{l=0}^{M_0} a_l^{(j')} V_{n,l}^{(j,j')}\,,
\end{equation}
with
\begin{equation}
	V_{n,l}^{(j,j')} = \frac{1}{2Q}\sum_{I(m,j)=j'} \left( \frac{|c_{m,j} z_m + d_{m,j}|}{\left(c_{m,j} z_m + d_{m,j}\right)}\right)^k (y^*_{m,j})^{\frac{k}{2}} \exp(2\pi \iu (l z^*_{m,j} - n x_m))\,,
\end{equation}
where $\sum_{I(m,j)=j'}$ denotes the sum over all $0\leq m \leq 2Q-1$ for which $I(m,j) =j'$. We therefore get
\begin{equation}
	\sum_{j'=1}^{n(c)}\sum_{l=0}^{M_0} a^{(j')} \tilde{V}_{n,l}^{(j,j')} = 0 \,,
\end{equation}
where
\begin{equation}
	\tilde{V}_{n,l}^{(j,j')} = V_{n,l}^{(j,j')} - \delta_{j,j'}\delta_{n,l} Y^{\frac{k}{2}} \exp(-2\pi n Y)\,,
\end{equation}
which we can again solve by imposing a normalization on the expansion at the cusp at infinity.

\subsection{A block-factored formulation of Hejhal's method}
\label{sec:hejhal_factored}
The matrix $V$, whose entries are given by Eq.~\eqref{eq:V_nl}, can be written as the
matrix product of two matrices (see for example Voight and Willis
\cite{voight_willis} who used an analogous factorization for a similar problem)
\begin{equation}
	\label{eq:V_factored}
	V = J\cdot W\,,
\end{equation}
with
\begin{equation}
	\label{eq:J}
	J_{n,m} = \frac{1}{2Q}\left( \frac{|c_m z_m + d_m|}{\left(c_m z_m + d_m\right)}\right)^k \exp(-2\pi \iu n x_m)\,,
\end{equation}
and
\begin{equation}
	\label{eq:W}
	W_{m,l} = (y^*_m)^{\frac{k}{2}} \exp(2\pi \iu l z^*_m)\,.
\end{equation}
Analogously, we can write $\tilde{V}_{n,l}$ whose entries are given by
Eq.~\eqref{eq:V_tilde_nl} as
\begin{equation}
	\tilde{V} = J\cdot W - D\,,
\end{equation}
where $D$ is a diagonal matrix whose entries consist of $Y^{\frac{k}{2}} \exp(-2\pi n Y)$.
For subgroups with more than one cusp, $V$ can be factored into a \emph{block-factored} form. For example for two cusps, we would get a matrix of the form
\begin{equation}
	\label{eq:V_tilde_block_factored}
	\tilde{V} =
	\begin{pmatrix}
		J^{(1,1)}\cdot W^{(1,1)} & J^{(1,2)}\cdot W^{(1,2)}\\
		J^{(2,1)}\cdot W^{(2,1)} & J^{(2,2)}\cdot W^{(2,2)}
	\end{pmatrix}
	-
	\begin{pmatrix}
		D^{(1)} & 0\\
		0 & D^{(2)}
	\end{pmatrix}\,.
\end{equation}
The same approach works analogously for more than two cusps. The factorization of the
involved matrices not only simplifies the expressions but can also significantly
improve the performance as we will discuss in the next section.

\section{Numerical computation of modular forms}
\label{sec:noncong_numerics}
We now discuss how Hejhal's method can be applied to compute numerical approximations of Fourier coefficients of modular forms on noncongruence subgroups. Because the matrices $J$ and $W$ can be efficiently constructed (for example by computing the corresponding powers through recursive multiplications), the computational bottleneck of Hejhal's method when working with a $q$-expansion basis is given by the linear algebra involved in the construction of $V$ and the linear solving. For this reason we survey different approaches for this task and present a new iterative mixed-precision approach that speeds up the linear solving significantly.

\begin{remark}[Implementational details]
	The algorithms discussed in this section have been implemented as a \textsc{Sage} \cite{sage} program. To compute the pullbacks we made use of the routines available in \textsc{Psage} \cite{psage}. For the LLL algorithm we used the implementation of \textsc{Pari} \cite{pari}. We also used \textsc{NumPy} \cite{numpy} and \textsc{SciPy} \cite{scipy} for double-precision computations. The arbitrary precision arithmetic has been performed using \textsc{Arb} \cite{arb} which is particularly useful in our application because of its highly optimized linear algebra routines \cite{arb_linear_algebra}.

	We plan to make our implementations publicly accessible in the future by contributing them to the \textsc{Psage} library.
\end{remark}

\subsection{The classical approach}
\label{sec:classical_solving}
Hejhal \cite{hejhal1999} and the majority of previous works
constructed the matrix~$\tilde{V}$ explicitly by performing
$\bigO(N^3)$ matrix multiplications between $J$ and $W$ and afterwards
used a $\bigO(N^3)$ direct solving technique to solve the resulting
linear system of equations. It is out of the question to compute
larger examples using arbitrary precision arithmetic with this
approach.

\subsection{The non-preconditioned Krylov approach}
\label{sec:non_precond_gmres}
To overcome the $\bigO(N^3)$ construction of the matrix $\tilde{V}$, Klug-Musty-Schiavone-Voight
\cite{klug_musty_schiavone_voight_2014} used a Krylov solving technique which only
requires the computation of matrix-vector-products, which means that $\tilde{V}$ can
be left in a block-factored form (we remark that Krylov solving techniques are also applied in the numerical method of \cite{monien_j2,monien_co3}). The convergence rate of this iterative solving
technique can be improved by scaling each column of $\tilde{V}$ by the diagonal term
which clusters the eigenvalues closer together. This gives (recall that
right-multiplying a matrix by a diagonal matrix corresponds to scaling its columns by
the diagonal entries)
\begin{align}
	\tilde{V}_\textrm{sc}
	:&= \tilde{V}
	\cdot
	{\begin{pmatrix}
			D^{(1)} & 0\\
			0 & D^{(2)}
	\end{pmatrix}}^{-1}\,,\\
	&= \left( \begin{pmatrix}
		J^{(1,1)}\cdot W^{(1,1)} & J^{(1,2)}\cdot W^{(1,2)}\\
		J^{(2,1)}\cdot W^{(2,1)} & J^{(2,2)}\cdot W^{(2,2)}
	\end{pmatrix}
	-
	\begin{pmatrix}
		D^{(1)} & 0\\
		0 & D^{(2)}
	\end{pmatrix}
	\right) \cdot
	{\begin{pmatrix}
			D^{(1)} & 0\\
			0 & D^{(2)}
	\end{pmatrix}}^{-1}\,,\\
	\label{eq:V_tilde_sc_factored}
	&=
	\begin{pmatrix}
		J^{(1,1)}\cdot W^{(1,1)} & J^{(1,2)}\cdot W^{(1,2)}\\
		J^{(2,1)}\cdot W^{(2,1)} & J^{(2,2)}\cdot W^{(2,2)}
	\end{pmatrix}
	\cdot
	{\begin{pmatrix}
			D^{(1)} & 0\\
			0 & D^{(2)}
	\end{pmatrix}}^{-1}
	-
	\begin{pmatrix}
		\mathbb{1} & 0\\
		0 & \mathbb{1}
	\end{pmatrix}\,.
\end{align}
The linear system therefore becomes
\begin{equation}
	\underbrace{\tilde{V} \cdot D^{-1}}_{=\tilde{V}_\textrm{sc}} \underbrace{\cdot D \cdot c}_{:=c'} = b\,,
\end{equation}
which we can solve for $c'$ to compute $c = D^{-1}c'$.

This approach typically runs faster compared to the classical approach. Its
limitation is however that the iteration count (i.e., the number of iterations until
convergence has been achieved) can become very high for involved problems with large
dimensions of $\tilde{V}$.

\subsection{The mixed precision iterative approach}
\label{sec:mixed_prec_solving}
To reduce the iteration count of an iterative method one typically attempts to find a preconditioner matrix $M$ to instead solve the linear system of equations
\begin{equation}
	M\cdot \tilde{V}_\textrm{sc}\cdot c' = M\cdot b\,,
\end{equation}
with improved convergence rate. However obtaining such a preconditioner appears
non-trivial for our application because $\tilde{V}_\textrm{sc}$ is non-hermitian,
non-symmetric and dense. In fact, we do not even know $\tilde{V}$ explicitly and, as
discussed before, constructing it is a $\mathcal{O}(N^3)$ operation so we would be in
the same order of magnitude as just applying a direct method to compute the
solution. The key observation to resolve this dilemma is that $\tilde{V}_\textrm{sc}$
can be safely inverted at a low precision. This can be seen from
Eq.~\eqref{eq:V_tilde_sc_factored}: The entries of the block matrices $W$ decay and
become effectively zero from a low precision perspective. Because $J$ does not change
the order of magnitudes, $J\cdot W$ also has decaying columns. However, by
subtracting the unit diagonal matrix we ensure that each column has at least one
entry that is non-zero. This means that if the Fourier expansion order $M_0$ is taken
to be very large, we asymptotically approach the unit matrix which is (and remains)
well-conditioned for inversion. We can therefore set the preconditioner $M$ to a
low-precision inverse (or something similar) of $\tilde{V}_\textrm{sc}$.

Such an approach uses mixed-precision arithmetic which is a relatively
recent concept that arose in HPC (high performance computing). For an
overview of different methods and applications utilizing
mixed-precision arithmetic we refer to
\cite{mixed_precision_survey}. The basic concept of mixed-precision
arithmetic is to perform computationally expensive parts of an
algorithm in faster low-precision arithmetic without sacrificing
precision of the end result. In the context of iterative solvers, it
has been shown and analyzed that low-precision inverses (of
potentially even highly ill-conditioned matrices) can serve as good
preconditioners for iterative methods \cite{rump_approximate_inverses,
  rump_approximate_inverses2, gmres_ir}. (In general, inverses are
good preconditioners because one approximates the unit matrix which
has the maximal clustered eigenvalue spectrum. If one knows the
inverse to full precision then the problem can obviously be solved in
one iteration but obtaining such an inverse is more expensive and
numerically unstable than solving the problem directly.) So far,
applications of mixed-precision arithmetic have typically replaced
double (64-bit) arithmetic with 32-/16-bit arithmetic which has faster
memory bandwidth and vectorization potential and is supported by
specialized hardware such as GPUs and Tensor Cores. The difference in
performance when switching between \emph{hardware supported} types
such as double-precision and arbitrary precision arithmetic is even
bigger: we find that reproducing a computation performed at
double-precision with the same precision using arbitrary precision
arithmetic takes about three orders of magnitude longer. Our approach
is therefore to construct $\tilde{V}$ explicitly in 64-bit
double-arithmetic and to compute an approximate inverse using a direct
method. Because these operations are performed in double-arithmetic
their contribution to CPU-time can be neglected in our examples.

\subsubsection{Preconditioned GMRES}
\label{sec:preconditioned_gmres}
\begin{figure}
	\centering
	\includegraphics[width=\columnwidth]{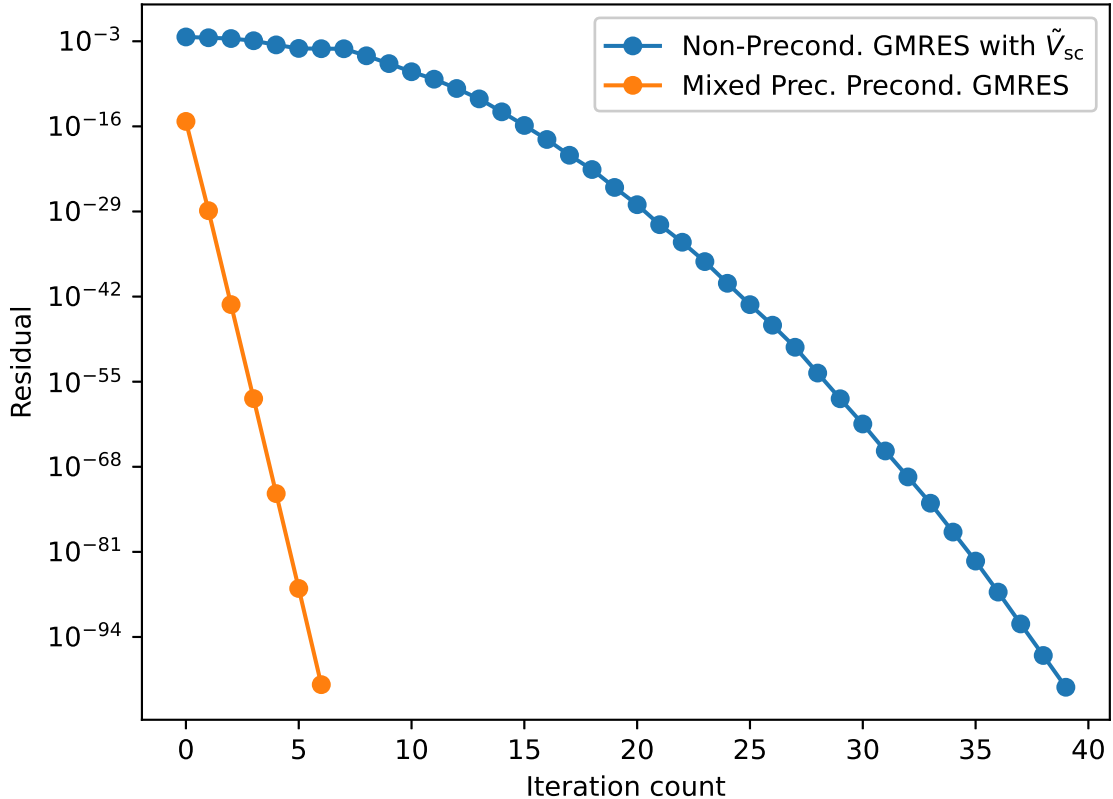}
	\caption{Comparison of precond. and non-precond. GMRES for the application of computing $f_0 \in S_4(G)$ where $G$ is a noncongruence subgroup with signature $(16, 0, 4, 0, 1)$ generated by $\sigma_S= (1\,10)(2\,14)(3\,7)(4\,12)(5\,16)(6\,8)(9\,15)(11\,13)$ and $\sigma_R= (1\,11\,14)(2\,15\,10)(3\,8\,7)(4\,13\,12)(5)(6\,9\,16)$ to 100 digits precision (taking $M_0 = 533$). The precond. version reduces the iteration count from 40 to 7 iterations.}
	\label{fig:gmres}
\end{figure}

As an example of an iterative Krylov subspace solver we implemented GMRES \cite{gmres}. To precondition GMRES with a low precision inverse, we first construct $\tilde{V}_\textrm{sc}$ in double-precision (which we will denote by $\tilde{V}_\textrm{sc,double}$) and compute its LU-decomposition
\begin{equation}
	\bar{L}\cdot \bar{U} = \tilde{V}_\textrm{sc,double}\,,
\end{equation}
where $\bar{L}$ and $\bar{U}$ denote the $L$ and $U$ factors up to double-precision. To compute the action of the inverse of $\tilde{V}_\textrm{sc,double}$ it is beneficial not to form $\tilde{V}_\textrm{sc,double}^{-1}$ explicitly which is computationally expensive, ill-conditioned and destroys potential sparseness. A better approach is to use \cite{gmres_ir}
\begin{equation}
	\tilde{V}_\textrm{sc,double}^{-1} x = \bar{U}^{-1} \bar{L}^{-1} x\,.
\end{equation}
The actions of ${\bar{L}}^{-1}$ and ${\bar{U}}^{-1}$ on a vector can be computed using $\bigO(N^2)$ triangular solves. Although the inverse is never explicitly formed, we will for simplicity still refer to this approach as \emph{computing the inverse}. The algorithm for preconditioned GMRES is illustrated in Algorithm \ref{alg:preconditioned_gmres}.
\begin{algorithm}
	\caption{Algorithm for computing Fourier expansion coefficients using GMRES}
	\label{alg:preconditioned_gmres}
	\begin{algorithmic}[1] 
		\STATE Compute block-factored form of $\tilde{V}_\textrm{sc}$ at full precision\;
		\STATE Construct $\tilde{V}_\textrm{sc,double}$ at double-precision\;
		\STATE Compute $\bar{L}\cdot \bar{U} = \tilde{V}_\textrm{sc,double}$ at double-precision\;
		\STATE Cast $\bar{L}$, $\bar{U}$ to full precision\;
		\STATE Solve $\left(\bar{L}\cdot  \bar{U}\right)^{-1} \tilde{V}_\textrm{sc} \cdot a' = \left(\bar{L}\cdot  \bar{U}\right)^{-1} b$ at full precision using GMRES\;
		\STATE Return $a = D^{-1}\cdot a'$ at full precision
	\end{algorithmic}
\end{algorithm}
The benefit of this algorithm is that GMRES gains \emph{at least} 16 digits (assuming $\tilde{V}_\textrm{sc,double}$ is well-conditioned) during each iteration. The reason for this upper bound on the iteration count (at least heuristically) comes from the fact that the inverse is known to 16 digits precision which means that the solution can be refined to 16 digits precision during each iteration. This convergence rate is not only very fast but it is also remarkable that the upper bound on the iteration count is (in principle) independent on the problem and the size of the matrices involved. (We only say \emph{in principle} because we assume here that the inverse of the matrix can be computed to 16 digits precision.) An illustration of the convergence rates for preconditioned and non-preconditioned GMRES can be found in Fig. \ref{fig:gmres}.

\subsubsection{Mixed precision iterative refinement}
\label{sec:IR}
Because GMRES needs to form a Krylov subspace, the action of $\tilde{V}_\textrm{sc}$ on a vector needs to be evaluated at the target precision during each iteration which is (comparatively) expensive. An alternative iterative algorithm which does not create a Krylov subspace is given by \emph{iterative refinement}. Iterative refinement (IR) is a relatively old technique that has first been applied by Wilkinson \cite{wilkinson} in 1948 and can be viewed as Newton's method on the function $r(x) = A\cdot x -b$ \cite{demmel}. In our application, the low precision inverse can be used to iteratively refine the solution vector during each iteration. Because we do not form a Krylov subspace we can gradually increase the precision during each iteration and do not need to perform all iterations at full precision. We therefore not only switch between double-precision and arbitrary precision arithmetic but also select different bit precisions when using arbitrary precision arithmetic. This approach therefore makes even more use of \emph{mixed precision} and is highlighted in Alg. \ref{alg:iterative_refinement}.
\begin{algorithm}
	\begin{algorithmic}[1] 
		\STATE Compute block-factored form of $\tilde{V}_\textrm{sc}$ at full precision\;
		\STATE Construct $\tilde{V}_\textrm{sc,double}$ at double-precision\;
		\STATE Compute $\bar{L}\cdot \bar{U} = \tilde{V}_\textrm{sc,double}$ at double-precision\;
		\STATE Use $\bar{L}\cdot \bar{U}$ to solve $\tilde{V}_\textrm{sc} \cdot a' = b$ at 64-bit\;
		\FOR{i=0:\textrm{max\_iter}-1}
		\STATE Compute $r = b - \tilde{V}_\textrm{sc}\cdot a'$ at $(i+2)\cdot 16$ digits precision\;
		\STATE Use $\bar{L}\cdot \bar{U}$ to solve $\tilde{V}_\textrm{sc} \cdot d = r$ at 64-bit\;
		\STATE Compute $a' = a' + d$\;
		\IF{\textrm{converged}}
		\STATE break\;
		\ENDIF
		\ENDFOR
		\STATE Return $a = D^{-1}\cdot a'$ at full precision
	\end{algorithmic}
	\caption{Algorithm for computing Fourier expansion coefficients using mixed-precision iterative refinement}
	\label{alg:iterative_refinement}
\end{algorithm}

Working at lower precisions during the iterations offers performance benefits for two reasons: First, arbitrary precision arithmetic has asymptotic complexity w.r.t. the precision $p$ given by $\bigO(p\log(p)\log(\log(p)))$ \cite[Section 2.3]{brent_zimmermann_2010} which means that working at a lower precision improves the performance of the ring operations. Second, because the matrix $W$ has decaying columns, many terms can be neglected when performing the matrix-vector product at a low precision.

If the approximate inverse is computed to 16 digits precision then iterative refinement gains 16 digits precision during each iteration. Contrary to GMRES, the convergence rate can only be linear which means that the iteration count of IR is larger than or equal to the one of GMRES.

\subsubsection{Precond. GMRES vs. mixed precision iterative refinement}
As discussed in the previous sections, GMRES can have a lower iteration count than IR while the iterations of IR are on average \emph{cheaper} because they do not have to be performed at the target precision. It is therefore interesting to examine which of these tradeoffs is beneficial in practice. For the examples that we have considered we find that IR typically has lower running times because the superlinear convergence of GMRES often only becomes noticeable during the last iterations (especially for larger index examples). For this reason, we have used mixed-precision IR as the numerical solver throughout this work.

\subsubsection{Optimizing the action of $W$}
\label{sec:W_action}
The action of $W$ (given by Eq. \eqref{eq:W}) can be interpreted as the evaluation of a polynomial at different points $q_m^*$ times factors $(y^*_m)^{\frac{k}{2}}$. It is a well known result that evaluating a polynomial at different points can be achieved in $\mathcal{O}(N\cdot \ln(N)^2)$ asymptotic complexity (see for example \cite{vandermonde_multiplication}). This asymptotic growth comes however with a large constant which makes this algorithm slower in practice than the classical $\mathcal{O}(N^2)$ algorithms for the problems that are considered in this work (additionally, these asymptotically fast algorithms are usually quite ill-conditioned).

For the classical $\mathcal{O}(N^2)$ algorithms, the most common choice would be
Horner's method which evaluates a polynomial at a single point using $N$
multiplications and $N+1$ additions as well as $\mathcal{O}(1)$ storage space.
However, because the powers of $q_m^*$ decay relatively fast, it is in practice
significantly faster to make use of \textsc{Arb}'s optimized dot-product routines
\cite{arb_linear_algebra} which, among other technical optimizations, evaluate each
term at the lowest possible precision (note that smaller terms can be evaluated at a
lower precision than larger terms without effecting the precision of the result).
Additionally, the dot-product routines neglect all terms that do not affect the
result. This is particularly useful because the iterative refinement algorithm (see
Algorithm \ref{alg:iterative_refinement}) starts with significantly lower precisions
(starting from 32 digits) than the target precision which means that the polynomials
can on average be truncated to lower order with many terms being neglected. Recall
also that $M_0$ is chosen based on the lowest point inside the fundamental domain, so
quite pessimistically, which means that the polynomials for many $\tau_m^*$ converge
faster. We note however that the naive approach of applying the dot-product, which
assembles the entries of matrix $W$ and computes its action by using the dot-product
row-wise, is not ideal for two reasons: First, the construction of W is comparatively
expensive because it requires $N^2$ multiplications at full precision which cannot be
further sped up. Second, and more importantly, storing $W$ as a matrix requires $N^2$
storage space in memory which becomes inconvenient for larger problems. For this
reason, we use modular splitting (see for example \cite[Section
4.4.3]{brent_zimmermann_2010}) for which only some of the powers of $q_m^*$ need to
be precomputed and stored. Modular splitting evaluates a polynomial $P(x)$ by using
the relations
\begin{equation}
	P(x) = \sum_{n=0}^{N} a_n x^n = \sum_{l=0}^{j-1} x^l P_l(x) = \sum_{l=0}^{j-1} x^l \left(\sum_{m=0}^{k-1} a_{jm+l}y^m \right)\,,
\end{equation}
where $y = x^j$. By choosing $j$ and $k$ to be of size $\mathcal{O}(\sqrt{N})$, we hence only need to store $\mathcal{O}(N^{3/2})$ values and can evaluate $P_l(x)$ using dot-products. We remark that we do not use classical rectangular splitting here because we do not want the terms of $P_l(x)$ to be uniformly distributed in order to make best use of the dot-product optimizations. We find that using \textsc{Arb}'s dot product often leads to a speedup that is close to an order of magnitude compared to a naive Horner scheme.

\subsubsection{Optimizing the action of $J$}
\label{sec:J_action}
It is immediate to see that the entries of $J$ (given by Eq. \eqref{eq:J}) are uniform and cannot be truncated when working at a lower precision which makes matrix-vector multiplication of $J$ very slow compared to $W$. We note however that $J$ can be further factored into:
\begin{equation}
	\label{eq:J_factored}
	J = D_L \cdot F \cdot D_R\,,
\end{equation}
where
\begin{equation}
	(D_L)_{n',m} = \exp\left(\frac{\pi i (2Q-1)}{2Q}\cdot n'\right)\,,
\end{equation}
\begin{equation}
	F_{n',m} = \exp\left(\frac{-2\pi i}{2Q}\cdot n'\cdot m\right)\,,
\end{equation}
\begin{equation}
	(D_R)_{n',m} = \frac{1}{2Q}\left( \frac{|c_m z_m + d_m|}{\left(c_m z_m + d_m\right)}\right)^k \exp\left(\frac{\pi i M_s (2Q-1)}{2Q}\right) \exp\left(\frac{-2\pi i M_s}{2Q}\cdot m\right)\,.
\end{equation}
Here $M_s$ denotes the index of the first coefficient that is non-zero (in general,
$M_s$ depends on the cusp, so we would instead need to write $M_s(j)$, but for the
sake of simpler notation, we assume $M_s$ to be equal for all cusps here) and $n' :=
n-M_s$ with the property $0 \leq n' \leq M-M_s$. $D_L$ and $D_R$ are diagonal
matrices whose action can be computed in $\mathcal{O}(N)$ operations. The matrix~$F$
is similar to the matrix of the classical discrete Fourier transform (DFT), but with
(in general) some missing rows and columns. Nevertheless, we can compute the action
of $F$ through a DFT. To illustrate this, assume that $M=3$, $2Q=4$ (obviously, in
practice we require $Q > M$) and that we have a missing column at $m=2$. Then the
action of $F$ on a vector can be written as:
\begin{equation}
	\begin{pmatrix}
		1 & 1 & 1\\
		1 & \left(\zeta_4\right)^{-1} & \left(\zeta_4\right)^{-3}\\
		1 & \left(\zeta_4\right)^{-2} & \left(\zeta_4\right)^{-6}
	\end{pmatrix}
	\cdot
	\begin{pmatrix}
		x_0\\
		x_1\\
		x_2
	\end{pmatrix}\,,
\end{equation}
where $\zeta_4 = \exp\left(\frac{2\pi i}{4}\right)$ is the 4-th root of unity. This is equivalent to computing:
\begin{equation}
	\begin{pmatrix}
		1 & 1 & 1 & 1\\
		1 & \left(\zeta_4\right)^{-1} & \left(\zeta_4\right)^{-2} & \left(\zeta_4\right)^{-3}\\
		1 & \left(\zeta_4\right)^{-2} & \left(\zeta_4\right)^{-4} & \left(\zeta_4\right)^{-6}\\
		1 & \left(\zeta_4\right)^{-3} & \left(\zeta_4\right)^{-6} & \left(\zeta_4\right)^{-9}
	\end{pmatrix}
	\cdot
	\begin{pmatrix}
		x_0\\
		x_1\\
		0\\
		x_2
	\end{pmatrix}\,,
\end{equation}
and selecting the first $3$ entries of the output vector. Our strategy for computing
the action of $F$ on a vector is therefore to zero-pad all entries of the input
vector for which $I(m,j) \neq i$, perform a DFT and afterwards select the first $M$
entries of the output vector. The advantage of using a DFT for computing the action
of $F$ is that fast Fourier transform (FFT) algorithms are available which have
asymptotic complexity $\mathcal{O}(N \ln (N))$ \cite{cooley_tukey}. Contrary to the
polynomial multipoint evaluation algorithms that were mentioned in Section
\ref{sec:W_action}, the FFT algorithms typically only have a small asymptotic
constant. In practice, we found the running time to be approximately $c\cdot Q\ln
(Q)$ where $c < 10$ if the largest prime factor of $Q$ is reasonably small (we used
the implementation provided by \textsc{Arb} to compute the FFT which has been contributed by
Pascal Molin). Because we have a free choice of $Q>M$, we choose $Q$ to be
\emph{slightly} larger than $M$ and with small prime factors to speed up the FFT.
Comparing this to the direct approach of computing the action of $J$ through
matrix-vector multiplications which has complexity $\mathcal{O}(Q\cdot M_0)$ (the
exact operation count also depends on the number of cusps) it is typically much
faster to use FFTs and the bottleneck of the algorithm becomes the action of $W$.
Additional advantages of factoring $J$ into the form of Eq.~\eqref{eq:J_factored} are
that the memory consumption becomes much lower since we only need to store the
diagonals and $2Q$ roots of unity, and that we avoid the $N^2$ operation to compute
the entries of $J$.

\subsubsection{Construction of $\tilde{V}_\textrm{sc,double}$}
To construct $\tilde{V}_\textrm{sc,double}$, we truncate the columns of $W$ so that terms that are effectively zero at double-precision are ignored. Afterwards, we compute the action of $J$ on the remaining columns of $W$ through FFTs (using \textsc{NumPy} \cite{numpy}), similarly to Section \ref{sec:J_action}. The construction of $\tilde{V}_\textrm{sc,double}$ therefore requires $\mathcal{O}(N^2 \ln (N))$ operations at double-precision.

\subsubsection{Computing the LU-decomposition of $\tilde{V}_\textrm{sc,double}$}
The matrix~$\tilde{V}_\textrm{sc,double}$ is sparse, since all of its entries that
are below the double machine epsilon are neglected. To compute its LU-decomposition
we therefore make use of the sparse linear algebra routines of \textsc{SciPy} \cite{scipy}. We
are unaware of the computational complexities of these routines (these should depend
on sparseness and structure of $\tilde{V}_\textrm{sc,double}$ and its LU factors) but
in practice they only account for negligible CPU time.

\subsubsection{Performing the LU-solves}
As discussed in the previous section, we use the sparse linear algebra routines of \textsc{SciPy} to compute a LU-decomposition of $\tilde{V}_\textrm{sc,double}$ in double arithmetic. When using this precomputed LU-decomposition to perform the solves inside the iterative refinement algorithm one needs to be careful not to over-/underflow the double exponent range which is finite and can be easily exceeded for elements inside the residue vectors. One way to avoid underflows is to convert the LU-decomposition to 53-bit \textsc{Arb}s which have unlimited exponent range. Storing the LU-factored matrix as an \textsc{Arb}-matrix is however quite memory consuming because \textsc{Arb} currently does not offer sparse matrices and because the memory footprint of a \textsc{Arb} object is higher than that of a double. A preferable approach is based on the observation that the input vectors for the LU-solves have relatively uniformly distributed entries. For this reason, we scale all entries by a constant factor $2^e$ to put them inside the double-range, convert them to doubles, perform the LU-solve in double arithmetic using \textsc{SciPy}, convert the result back to \textsc{Arb} and scale the result back. This approach uses significantly less memory and is faster.

\subsubsection{Restarting the algorithm}
Because the iterative refinement algorithm does not need to form a Krylov-subspace,
it can be restarted without losing convergence. One approach that we have
experimented with gradually increases the values of $Q$ and $M_0$. For example if a
target precision of 500 digits is to be reached, one can first choose $Q$ and $M_0$
so that convergence is reached up to 100 digits precision. One can then afterwards
use these approximations of the lower coefficients up to 100 digits precision to
restart the algorithm with a larger choice of $Q$ and $M_0$ to refine the residue
from $10^{-100}$ to $10^{-250}$ and afterwards again to go from $10^{-250}$ to
$10^{-500}$. The performance that one can gain from this approach however seems to be
quite limited because the bottlenecks are given by the last iterations anyways.
Additionally, each restart creates some extra computations to set up $J$, $W$ and the
preconditioner. Although some restarting configurations exist that are faster than
simply starting with the target values of $Q$ and $M_0$, the performance impact is
very minor and finding these configurations can be inconvenient which is why we have
not applied this approach for our computations.

\subsubsection{Performance comparison to previous methods}
To examine how the different approaches perform in practice we ran several benchmarks that numerically compute a modular forms on a noncongruence subgroups at different precisions. The results of these benchmarks can be found in Tab. \ref{tab:higher_index_benchmarks}. We report the CPU times and peak memory usages of the program. All implementations are highly optimized from a technical perspective. The \emph{classical} version
follows the approach of Section \ref{sec:classical_solving} (we used \textsc{Arb}'s implementations for the matrix multiplications and LU decompositions). The \emph{non-precond. GMRES} version follows the approach of Section \ref{sec:non_precond_gmres} with GMRES as a Krylov solver. The \emph{mixed precision IR} version uses the mixed precision iterative refinement approach with optimized actions of $J$ and $W$ that was presented in Section \ref{sec:mixed_prec_solving}. The benchmarks where taken on a \texttt{Intel Xeon E5-2680 v4 @ 2.40GHz} CPU and ran on a single thread.
\begin{table}
	\centering
	\begin{tabular}{c|c|c|c}
		Digit precision ($M_0$) & Classical & Non-precond. GMRES & Mixed Precision IR\\
		\hline
		100  (533) & 7min15s (1.91GB) & 1min38s (1.5GB) & \textbf{7s (0.32GB)}\\
		200 (1043) & 1h10min (9.48GB) & 15min3s (5.84Gb) & \textbf{43s (0.47GB)}\\
		400 (2061) & - & 4h25min (29.19GB) & \textbf{6min19s (0.94GB)}
	\end{tabular}
	\caption {Benchmarks for the numerical computation of $f_0 \in S_4(G)$ where $G$ is a subgroup of signature $(17, 0, 3, 1, 2)$ that is generated by $\sigma_S = (1)(2\,4)(3\,7)(5\,10)(6\,11)(8\,14)(9\,15)(12\,13)(16\,17)$ and $\sigma_R = (1\,7\,4)(2\,11\,10)(3\,15\,14)(5)(6\,12\,13)(8\,17\,9)(16)$.}
	\label{tab:higher_index_benchmarks}
\end{table}
As one can see, the mixed-precision algorithm outperforms the other algorithms in all
categories and runs more than $40$ times faster than non-precond.\ GMRES at 400
digits precision while consuming significantly less memory. For larger examples this
ratio becomes even bigger because the IR approach has a lower asymptotic complexity.

\subsubsection{Numerical stability for large examples}
\begin{figure}
	\centering
	\includegraphics[width=\columnwidth]{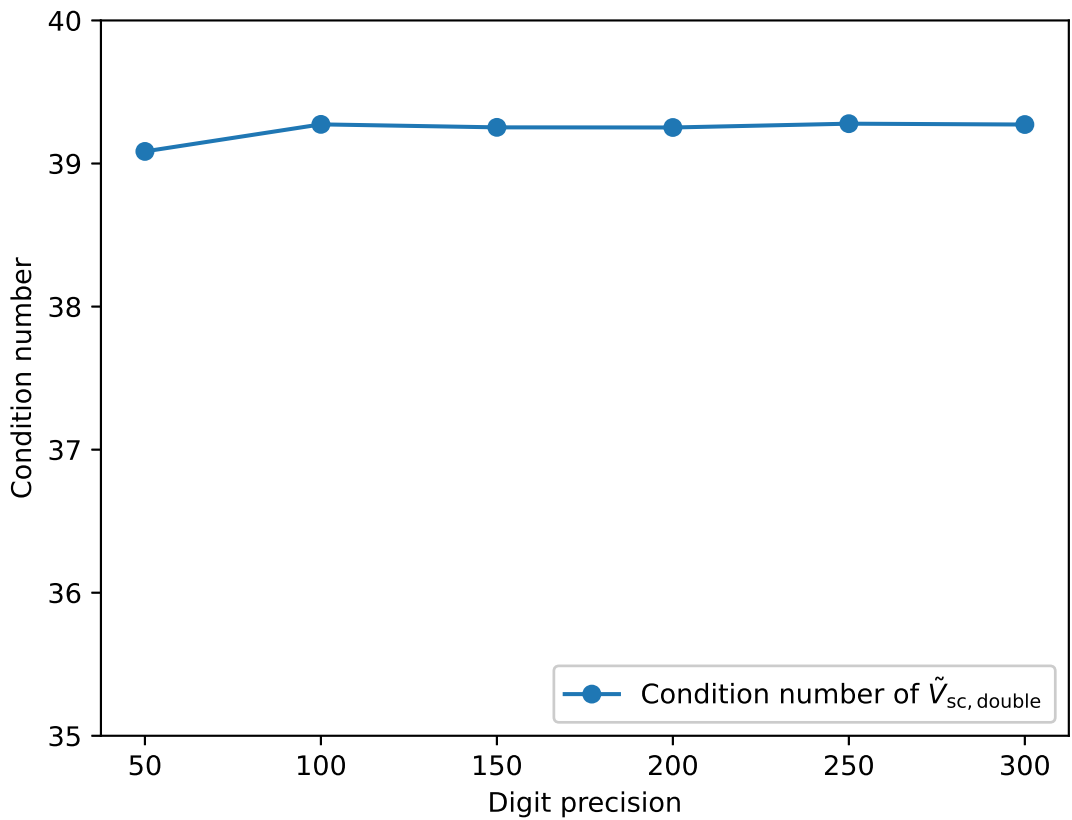}
	\caption{Illustration of the condition number of $\tilde{V}_\textrm{sc,double}$
	for varying target precisions. For the example we used the cusp form that was
	considered in Tab.~\ref{tab:higher_index_benchmarks}.}
	\label{fig:condition_numbers}
\end{figure}

Increasing the target precision (and following from that the values of $M_0$ and $Q$)
does not affect the condition number of $\tilde{V}_\textrm{sc,double}$ (up to some
noise), as illustrated in Fig.~\ref{fig:condition_numbers}. This seems to be caused
by the fact that the additionally added columns are similar to those of a unit-matrix.

\begin{figure}
	\centering
	\includegraphics[width=\columnwidth]{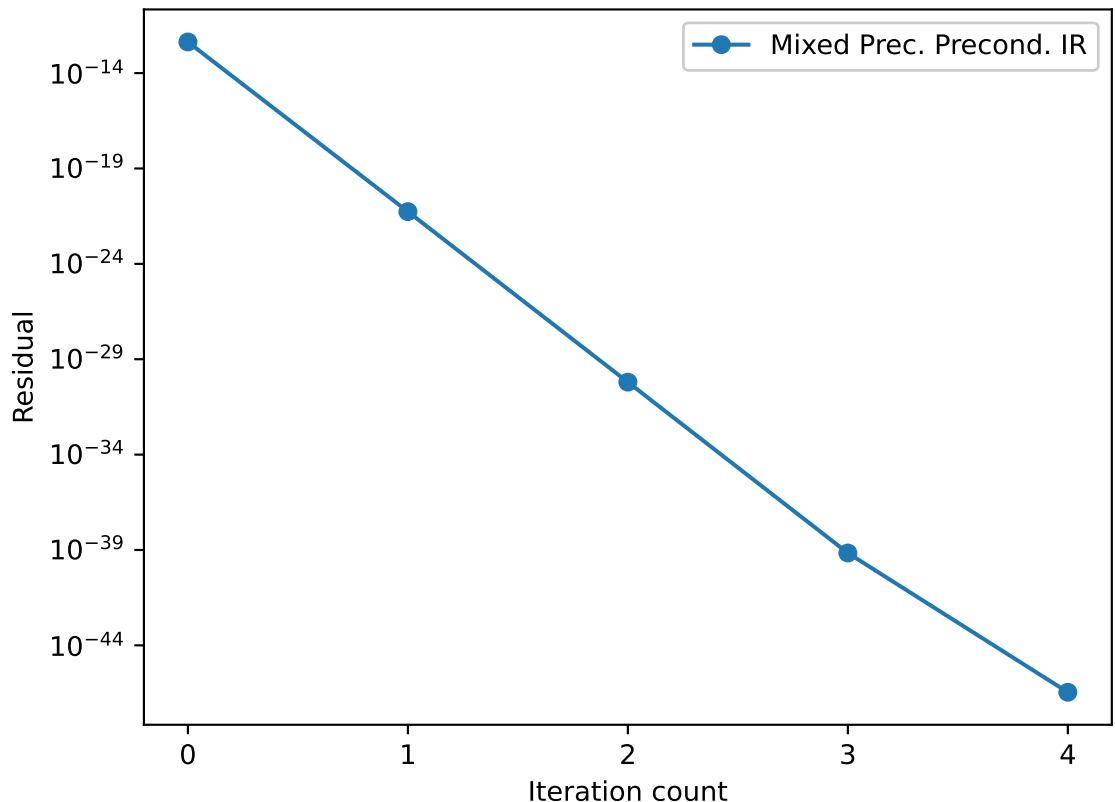}
	\caption{Illustration of the iterative computation of $f_0 \in S_2(\Gamma_0(120))$ to 50 digits precision using mixed precision IR (taking $M_0 = 2725$).}
	\label{fig:Gamma0_120}
\end{figure}

The index and number of cusps of the considered subgroup affect the conditioning more
noticeably. Although large index examples have not been the focus of this work it is
therefore interesting to examine if they are well-conditioned enough to apply mixed
precision iterative refinement on them as well. For this we consider the subgroup
$\Gamma_0(120)$ of signature $(288, 17, 16, 0, 0)$. (This is obviously a congruence
subgroup for which efficient non-numerical algorithms exist from which we can get
exact solutions. This makes it a useful example to test the numerical stability of
the algorithm. We also remark that the numerical method does not distinguish between
congruence and noncongruence subgroups which means that we can expect the same
results to hold for noncongruence subroups.) It is immediate to see that with an
index of 288 and 16 cusps of which the largest one has width 120, $\Gamma_0(120)$ is
significantly larger than the other considered examples. As a test of our algorithm
we performed the numerical computation of $f_0 \in S_2(\Gamma_0(120))$ to 50 digits
precision. To achieve convergence we take $M_0 = 2725$ which means that the resulting
linear system of equations is of dimension $43600\times 43600$ which is enormous in
the context of arbitrary precision arithmetic. Still, we found that iterative
refinement converges fast as can be seen in Fig. \ref{fig:Gamma0_120}. Contrary to
the other examples, the size of $\tilde{V}$ reduces the precision gain per iteration
to about 9 digits per iteration instead of 16. We also found that the resulting
coefficients have \emph{only} been computed to about 44 digits precision instead of
50 which can however obviously easily be overcome by setting a buffer for large index
examples. The computation used 60GB of memory and took 2h and 30min of CPU time. We
also remark that, contrary to the other computations, we had to use dense linear
algebra to perform the LU decomposition because the sparse routines returned a memory
error. We conclude that mixed-precision iterative refinement can be efficiently
applied to large index examples as well.

\subsubsection{Overall complexity of the algorithm}
Studying the complexity of the mixed-precision IR algorithm is relatively difficult.
First of all, it makes sense to ignore all computations that can be performed in
double arithmetic, because due to their technical optimization, their contribution
can be neglected compared to the parts that use arbitrary precision arithmetic (at
least for the scale of problems that are considered in this work and taking the limit
$N\rightarrow \infty$ would lead to conditioning problems at some point anyways).
When analyzing the performance with respect to $N$ (we use $N$ synonymously for $Q$
and $M_0$ because these are usually proportional to one another), the asymptotic
bottleneck both in theory and in practice is given by the action of $W$. The
complexity of this computation is $\mathcal{O}(N^2)$ (at least in practice, as
discussed in Section \ref{sec:W_action}, the theoretical asymptotic complexity is
$\mathcal{O}(N\cdot \ln(N)^2)$) but comes with a very small constant due to the
decaying columns of $W$ and the relatively small iteration count. We also remark
that, contrary to most iterative methods, the iteration count of our method only
depends on the precision and is hence truly independent of $N$.

A more meaningful quantity would be the bit-complexity of the algorithm. This does however seem impossible to calculate due to the constantly varying working precisions, floating point types and decay rates of the dot-product terms.

\subsection{Examples}
We now illustrate the approach of Section \ref{sec:mixed_prec_solving} for some examples.
\begin{example}
	\sloppy Let $G$ be the (randomly selected) noncongruence subgroup with signature $(16, 1, 2, 0, 1)$ that is generated by $\sigma_S = (1\,4)(2\,5)(3\,8)(6\,11)(7\,10)(9\,14)(12\,15)(13\,16)$ and $\sigma_R = (1)(2\,10\,11)(3\,7\,14)(4\,8\,5)(6\,16\,15)(9\,13\,12)$. Following from this, we get that $\sigma_T = (1\,8\,7\,11\,16\,12\,6\,2\,4)(3\,5\,10\,14\,13\,15\,9)$ which means that the cusp at infinity has width 9. Note that $\textrm{dim}(S_2(G)) = 1$. We use the approach of Section \ref{sec:mixed_prec_solving} to compute $f_0 \in S_2(G)$ to 150 digits precision. This computation takes about 5s on a standard CPU. By recognizing $a_2^9$ as an algebraic number using the LLL algorithm \cite{LLL} we find that $K = \QQ(v)$, where
	\begin{equation}
		v^3 - 6v - 16 = 0\,,
	\end{equation}
	with embedding $v = -1.647426...+1.463572...\iu$. Choosing
	\begin{equation}
		u = \left(\frac{3^3\cdot 49667\cdot 1452815993}{2^{45}\cdot 7^7\cdot 137^9} - \frac{3^4\cdot 5\cdot 14543 \cdot 393024407}{2^{47}\cdot 7^7\cdot 137^9}v  -\frac{3^4\cdot 167\cdot 9697\cdot 1862489}{2^{51}\cdot 7^7\cdot 137^9}v^2\right)^{1/9}\,,
	\end{equation}
	and applying the LLL algorithm again we recognize the Fourier expansion to be given by
	\begin{equation}
		f(q_9) = q_9 + (822u)q_9^2 + ((-68028v^2 - 253920v - 445797)u^2)q_9^3 + ...\,,
	\end{equation}
	up to the 23rd order (higher orders can be recognized by increasing the target precision). Although this result is based on very high heuristic evidence it is not yet formally proven (it should in principle be possible to prove the result through the curve but we do not carry this out here).
\end{example}
\begin{example}
  \sloppy A more complicated example is given by the subgroup with
  signature $(13, 1, 1, 1, 1)$ that is generated by
  $\sigma_S = (1)(2\,4)(3\,7)(5\,10)(6\,9)(8\,12)(11\,13)$,
  $\sigma_R= (1\,7\,4)(2\,9\,10)(3\,6\,12)(5\,8\,13)(11)$ and
  $\sigma_T= (1\,7\,6\,10\,8\,3\,4\,9\,12\,13\,11\,5\,2)$. Computing
  $f_0 \in S_2(G)$ to 1000\ digits precision which takes about 90\
  minutes of CPU time on a \texttt{Intel Xeon E5-2680 v4 @ 2.40GHz} we find that $K = \QQ(v)$, where
	\begin{equation}
		 v^{10} - 3v^9 + 5v^8 - 12v^7 + 24v^6 - 46v^5 + 68v^4 - 60v^3 + 96v^2 - 144v + 72 = 0\,,
	\end{equation}
	with embedding $v = 1.068141... + 0.135042...\iu$ and compute the corresponding cusp form up to the 25th order (the resulting expressions become too large to be displayed here). This example would not be feasible to compute with the previous methods.
\end{example}

\section{Computation of modular forms on genus zero subgroups}
\label{sec:genus_zero_belyi_maps_and_forms}

The methods of Section \ref{sec:noncong_numerics} can obviously be applied to
numerically compute modular forms on subgroups of arbitrary genus. In this section we
discuss a different approach that is restricted to subgroups of genus zero, for which
the field of modular functions is generated by a single function, called the
\emph{Hauptmodul} (see Section~\ref{sec:hauptmoduls}). The methods described in this
section can be used to obtain \emph{rigorous} results.

\subsection{Computing genus zero Belyi maps}
\begin{theorem}[Atkin-Swinnerton-Dyer]
	\label{th:asd}
	A necessary and sufficient condition that $f(\tau)$ is a modular function on a
	subgroup of finite index in $\Gamma$ is that $f(\tau)$ should be an algebraic
	function of $j$ and that its only branch points should be branch points of order
	2 at which $j=1728$ and branch points of order 3 at which $j=0$, and branch
	points at which $j$ is infinite.
	\begin{proof}
		See Atkin-Swinnerton-Dyer \cite[Theorem 1]{asd}.
	\end{proof}
\end{theorem}

In particular, note that $j$, when viewed as a function on the
modular curve $X(G)$ of some finte index subgroup $G \subgroup \Gamma$,
gives an example of a \emph{Belyi map}.
\begin{definition}[Belyi Map]
	Let $X$ be a compact Riemann surface. Then a holomorphic function
	\begin{equation}
		f\, : \, X \rightarrow \mathbb{P}^1(\CC)\,,
	\end{equation}
	is said to be a \emph{Belyi map} if it is unramified away from three points.
\end{definition}
Belyi maps inherit their name from a famous theorem by Belyi~\cite{belyi}
\begin{theorem}[Belyi]
	A compact Riemann surface $X$ (equivalently an algebraic curve) over~$\CC$ can be
	defined over $\overline{\QQ}$ if and only if there exists a Belyi map on $X$.
	\begin{proof}
		See Belyi \cite[Theorem 1]{belyi2}.
	\end{proof}
\end{theorem}
Belyi maps and their computation is an interesting subject on their own with numerous
applications in number theory and algebraic geometry, for an overview we refer to the
survey of Sijsling and Voight~\cite{sijsling_voight}.

Let $G$ be a finite index subgroup of $\Gamma$. Then the covering map
\begin{equation}
	R\, : \, X(G) \rightarrow X(\Gamma) \overset{j}{\cong}
	\mathbb{P}^1(\CC)\,,
\end{equation}
is a Belyi map, where $X(G)=G\backslash \Hbar$ is the modular curve. If $G$ is a genus zero subgroup then the covering map $R(j_G)$ is a rational function in $j_G$, and branches over the images of the elliptic points and cusps as well. This means that $R$ can be written as
\begin{equation}
	\label{eq:covering_map}
	R(j_G) = \frac{p_3(j_G)}{p_c(j_G)} = 1728 + \frac{p_2(j_G)}{p_c(j_G)}\,.
\end{equation}
The ramification structure (i.e., the roots of $p_2$, $p_3$ and $p_c$) can be determined from the cycle type of $\sigma_S$, $\sigma_R$ and $\sigma_T$. Let us illustrate this on an example.
\begin{example}[Determining ramification structure from permutation triple]
Let $G$ be a noncongruence subgroup with signature $(7, 0, 2, 1, 1)$ corresponding to
the permutations $\sigma_S = (1)(2\,4)(3\,5)(6\,7)$, $\sigma_R =
(1\,5\,4)(2\,7\,3)(6)$ and $\sigma_T = (1\,5\,2)(3\,4\,7\,6)$. By definition $p_2$
needs to be of the form
	\begin{equation}
		p_2(j_G(\tau)) = \prod_{i=1}^{7}\left(j_G(\tau)-j_G(e_{2,i})\right)\,,
	\end{equation}
	where we denote $e_{2,i}$ to be the elliptic point of order two, located at coset
	of index $i$. Because some of the values of $j_G(e_{2,i})$ at the elliptic points
	are equal, we can write this as
	\begin{equation}
		p_2(j_G(\tau)) = (j_G(\tau)-j_G(e_{2,1})) (j_G(\tau)-j_G(e_{2,2}))^2 (j_G(\tau)-j_G(e_{2,3}))^2 (j_G(\tau)-j_G(e_{2,6}))^2\,.
	\end{equation}
	This means that $p_2$ can be written in the form
	\begin{equation}
		p_2(j_G) = \left(j_G^3+A_2 j_G^2+ B_2 j_G + C_2\right)^2\left(j_G+D_2\right)\,,
	\end{equation}
	where (by Belyi's theorem) $A_2,B_2,C_2,D_2 \in \overline{\QQ}$. Analogously,
	$p_3$ and $p_c$ can be factored into
	\begin{equation}
		p_3(j_G) = \left(j_G^2 + A_3 j_G + B_3 \right)^3\left(j_G+C_3\right)\,,
	\end{equation}
	and
	\begin{equation}
		p_c(j_G) = \left(j_G+A_c\right)^4\,,
	\end{equation}
	where the roots are given by $j_G(e_{3,i})$, resp.\ $j_G(c_i)$.
\end{example}
Once the structure of $p_2$, $p_3$ and $p_c$ has been determined, we can transform Eq. \eqref{eq:covering_map} into
\begin{equation}
	\label{eq:belyi_poly_vanishing}
	P(j_G) := p_3(j_G) - p_2(j_G) - 1728 p_c(j_G) = 0\,,
\end{equation}
where $P(j_G)$ is a polynomial whose coefficients are defined over symbolic expressions.
The coefficients of $P(j_G)$ need to vanish which gives $\textrm{deg}(P) =
\textrm{deg}(p_3) = \textrm{deg}(p_2)$ polynomial equations in the unknowns
$A_2,A_3,\dots$. An additional equation is obtained by expanding $R(j_G)$ in
$j_G(q_N)$ and by asserting that the constant term is equal to 744 if the cusp width
at infinity is equal to one and vanishes otherwise.

One can attempt to solve these non-linear systems of equations directly, for example
by using Gröbner bases \cite[Section 2]{sijsling_voight}. This however quickly
becomes infeasible for all but the simplest examples. A much more efficient approach
is to use a numerical method to compute approximations of the evaluation of the
Hauptmodul at the elliptic points and the cusps. These
approximations can then be
used as starting values for Newton iterations to determine the unknown coefficients
to high precision. Afterwards the LLL algorithm can be applied to identify the
expressions as algebraic numbers. This approach has been suggested by
Atkin-Swinnerton-Dyer \cite{asd} and its effectiveness has been demonstrated by
the second author \cite{monien_j2,monien_co3} who used this approach to compute Belyi maps for
genus zero noncongruence subgroups of large index and degree of the number field.
Similar approaches that use approximations of modular forms as starting values for
Newton iterations have been used in
\cite{selander-stroembergsson,klug_musty_schiavone_voight_2014,belyi_db}.

\subsubsection{Obtaining starting values for Newton's method}
Throughout this work we have computed the starting values for Newton's method by
using the algorithm that is described in Section~\ref{sec:noncong_numerics}. The
Fourier expansion of the Hauptmodul at infinity can be normalized to be of the form
$q_N^{-1} + 0 + a_1q_N + a_2q_N^2 + ...$. The values of the Hauptmodul at the other
cusps are finite which means that its expansions are of the form $a_0 + a_1q_{N_c} +
...$.
\begin{remark}
	It is important to note that the $q_N^{-1}$-terms form the right-hand side of the
	linear system of equations and therefore do not enter~$\tilde{V}$. This means
	that the largest entries for each column of $\tilde{V}$ are still located on the
	diagonal and hence that the mixed-precision iterative solving techniques of
	Section~\ref{sec:noncong_numerics} can also be used to compute $j_G$.
\end{remark}
For the examples that have been considered in this work, it is sufficient to
numerically compute the Fourier expansion of the Hauptmodul to 50 digits of precision
(although computations at lower precision would have probably worked as well).
The evaluations of the Hauptmodul at the elliptic points can then be computed by
evaluating the Hauptmodul at $\gamma_i(\iu)$ and $\gamma_i(\rho)$ where~$\gamma_i$
denotes the coset representative of the corresponding coset (it is preferred to
choose the cusp expansion with the fastest convergence for the evaluation at these
points in order to maximize the precision). The values at the cusps outside infinity
are simply given by the constant terms in the cusp expansions.

\subsubsection{Applying Newton's method}
Once the starting values have been obtained the multivariate Newton method can be
used to improve the precision of these values. For simplicity we will use $x=j_G$ in
this section. We also use $[x^n]P$ to denote the coefficient of $x^n$ in $P$.
Then the Jacobian of the system of polynomial equations is given by a $(\mu+1)\times
(\mu+1)$ matrix (where $\mu$ denotes the index of $G$) that is of the form
\begin{equation}
	J(P) =
	\begin{pmatrix}
		\frac{\partial}{\partial A_2} [x^0]P & \frac{\partial}{\partial B_2} [x^0]P & \cdots\\
		\frac{\partial}{\partial A_2} [x^1]P & \frac{\partial}{\partial B_2} [x^1]P & \cdots\\
		\vdots & \vdots & \vdots\\
		\frac{\partial}{\partial A_2} [x^\mu]P & \frac{\partial}{\partial B_2} [x^\mu]P & \cdots
	\end{pmatrix}
	\,.
\end{equation}
Let $X^{[m]} \in \CC^{\mu}$ be the vector containing the numerical approximations of
the unknowns $A_2, B_2, \dots$ at the $m$-th iteration. Then we can use the update
steps
\begin{equation}
	X^{[m+1]} = X^{[m]} - [J(P(X^{[m]}))]^{-1}P(X^{[m]})\,,
\end{equation}
to iteratively increase the precision of the approximations of $X$. As is standard
with Newton's method, this procedure achieves quadratic convergence.

We remark that from a numerical perspective it is preferable to perform the update steps by solving the linear system of equations
\begin{equation}
	\label{eq:jacobian_solving}
	J(P(X^{[m]}))\cdot d^{[m]} = P(X^{[m]})\,,
\end{equation}
instead of computing the matrix inverse of the Jacobian (see the discussion in
Section \ref{sec:preconditioned_gmres}). Analogously to the iterative refinement, the
update steps are then given by
\begin{equation}
	X^{[m+1]} = X^{[m]} + d^{[m]}\,.
\end{equation}
We used \textsc{Arb}'s LU-decomposition to solve Eq.~\eqref{eq:jacobian_solving} (i.e., a
direct solving technique). For large index examples it might be preferable to perform
this solving iteratively, for example by using preconditioned GMRES (see
Section~\ref{sec:preconditioned_gmres}).

As an additional implementation detail we remark that instead of computing the
entries of the Jacobian matrix through symbolic computation of the partial
derivatives and their evaluations by plugging in the corresponding approximations of
the variables it is instead preferable to compute the columns of the Jacobian through
univariate polynomial multiplication which has been used by the second author
in~\cite{monien_j2,monien_co3}. To illustrate this, suppose that $P$ is of the form
\begin{equation}
	P = (a_0 + a_1 x + a_2 x^2 + ...)^{k_a} \cdot (b_0 + b_1 x + b_2 x^2 + ...)^{k_b}\cdot ... + ...\,,
\end{equation}
then
\begin{equation}
	\frac{\partial}{\partial a_i}P = k_a x^i (a_0 + a_1 x + a_2 x^2 + ...)^{k_a-1} \cdot (b_0 + b_1 x + b_2 x^2 + ...)^{k_b}\,.
\end{equation}
Constructing this polynomial by using multiplications of univariate polynomials over $\CC$ (for which we used \textsc{Arb}'s polynomial implementation) then yields in a polynomial whose coefficients correspond to a column of $J$, since
\begin{equation}
	\frac{\partial}{\partial a_i}[x^j]P = [x^j]\left(\frac{\partial}{\partial a_i}P\right)\,.
\end{equation}
By applying this procedure for all unknowns (and potentially reusing terms for optimization), all entries of $J$ can be assembled efficiently.

\subsubsection{Identifying the Belyi map}
Once the coefficients of the Belyi map have been computed to sufficient precision, the LLL algorithm can be used to identify $K$ and $u$.
\begin{example}
	Continuing the example of this section we find that the Belyi map is given by
	\begin{align}
	\begin{split}
		R(x) &= \frac{(x^2 + 444ux - 148284u^2)^3(x + 516u)}{(x + 462u)^4}\,,\\
		&= 1728 + \frac{(x - 996u)(x^3 + 1422ux^2 + 822204u^2x + 185029704u^3)^2}{(x + 462u)^4}\,,
	\end{split}
	\end{align}
	where $u = (2/823543)^{1/3}$ which means that $K = \QQ$.
\end{example}
We can verify that the result of the Belyi map is correct by confirming that
Eq.~\eqref{eq:belyi_poly_vanishing} holds for the recognized polynomials.

\subsection{Computing Fourier expansions of the Hauptmodul from the Belyi map}
\label{sec:hauptmodul_from_belyi}
The result of the Belyi map can be used to explicitly compute Fourier expansions of the Hauptmodul.

\subsubsection{Computing Fourier expansions at infinity}
\label{sec:hauptmodul_expansion_infinity}
We have seen that
\begin{equation}
	j = R(x) = \frac{p_3(x)}{p_c(x)}\,,
\end{equation}
which we hence need to solve for $x=j_G$. To do this we work with the reciprocal
\begin{equation}
	\frac{1}{j} = \frac{1}{R(x)} = \frac{p_c(x)}{p_3(x)} =: \frac{1}{R(1/\Tilde{x})}
	= \frac{p_c(1/\Tilde{x})}{p_3(1/\Tilde{x})}\,,
\end{equation}
where we set $\Tilde{x} := 1/x$. Expanding $1/R(1/\Tilde{x})$ as a power series in
$\Tilde{x}$ results in
\begin{equation}
	\sqrt[N]{\frac{1}{j}} = \sqrt[N]{\frac{1}{R(1/\Tilde{x})}} =: s(\tilde{x})\,,
\end{equation}
where $N$ denotes the width of the cusp at infinity and the roots denote the roots of
the power series. (We remark that in order to identify the correct embedding of the $N$-th root we compared the embeddings to the result of the numerical method of Section \ref{sec:noncong_numerics}.) The power series $s(\tilde{x})$ has valuation one and we can hence
compute the reversion
\begin{equation}
	\Tilde{x} = s^{-1}(\sqrt[N]{1/j})\,,
\end{equation}
to get
\begin{equation}
	x = 1/s^{-1}(\sqrt[N]{1/j})\,.
\end{equation}
Substituting the $q$-expansion of $\sqrt[N]{1/j}$ (which is a power series in $q_N$)
then gives the $q$-expansion of $j_G$ at
infinity.

\subsubsection{Computing Fourier expansions at other cusps}
\label{sec:hauptmodul_expansion_non_infinity}
To compute the Fourier expansion at a cusp $\ne i\infty$ we perform the transformation
\begin{equation}
	x \mapsto x + j_G(c_i) := \tilde{x}\,,
\end{equation}
where $j_G(c_i)$ denotes the evaluation at the cusp. Then
\begin{equation}
	\sqrt[N]{\frac{1}{j}} = \sqrt[N]{\frac{1}{R(\Tilde{x})}} =: s(\tilde{x})\,,
\end{equation}
where $N$ denotes the width of the considered cusp (not at infinity) and
\begin{equation}
	x = j_G(c_i) + s^{-1}(\sqrt[N]{1/j})\,.
\end{equation}

\subsubsection{Computing Fourier expansions over number fields}
To perform computations over number fields we introduce the number field $L = \QQ(w)$ over
which the coefficients of the Belyi map and the Fourier expansions are defined. If
the cusp width at infinity is equal to one then $K = L$. Otherwise we choose
$L$ to be the number field generated by $u$. The advantage of this choice of $L$ is that
one can efficiently convert its elements into \emph{u-v-factored} expressions (and
vice versa).

Once the Belyi map has been recognized explicitly over $L$, the expansions at
infinity can be computed by performing the arithmetic of Section
\ref{sec:hauptmodul_expansion_infinity} over $L$. For this we used the generic
routines provided by \textsc{Sage} \cite{sage}. The advantage of this approach is that the
Fourier coefficients of the Hauptmodul are rigorous. Note that expansions of cusps
outside infinity cannot in general be computed over $L$ because they are defined
over number fields $\QQ(v^{1/N_c})\QQ(w)$ where $N_c$ denotes
the cusp width of
the considered cusp outside infinity.

\subsubsection{Computing Fourier expansions over $\CC$}
\label{sec:hauptmodul_from_belyi_cc}
To compute Fourier coefficients of the Hauptmodul over $\CC$ (more precisely, using
arbitrary precision arithmetic) we use \textsc{Arb} \cite{arb} to perform the computations of
sections~\ref{sec:hauptmodul_expansion_infinity}
and~\ref{sec:hauptmodul_expansion_non_infinity}.
The bottleneck of these computations is the reversion of power series. We found that
series reversion in \textsc{Arb} is significantly faster than in \textsc{Sage} \cite{sage} or \textsc{Pari}
\cite{pari}. \textsc{Arb} has implemented the algorithm of
\cite{arb_series_reversion} which decreases the asymptotic complexity from
$\bigO(N^3)$ to $\bigO(N^{1/2}M(N)+N^2)$, where $M(N)$ denotes the complexity of
polynomial multiplication. \textsc{Arb} also provides implementations of the fast power series
composition algorithms of \cite{brent_kung_power_series_composition} which we use for
the substitutions.

We note however that the approach of sections \ref{sec:hauptmodul_expansion_infinity}
and \ref{sec:hauptmodul_expansion_non_infinity} can be very ill-conditioned which
means that one might have to use a higher working precision than the target
precision. This seems to be caused by the fact that the reversed series $s^{-1}$ can
have very large coefficients which makes the substitution ill-conditioned. We are
unaware of a transformation that improves the conditioning so the best we could come
up with is an approach where we choose the working precision \emph{sufficiently
large} in order to overcome the ill-conditioning. We first compute
$s^{-1}$ to low precision (typically to 64-bit, but not in double
precision because the exponents might over/underflow). The size of the resulting
coefficients gives an estimate of the required precision. If the computation fails at
the estimated precision, we attempt it again using a higher precision. The interval
arithmetic of \textsc{Arb} is very useful for this since it shows if the working precision
had been sufficiently high. While this strategy obviously always leads to correct
results, it is not very elegant and it would be useful to find a way to rewrite the
problem so that all computations can be done at the target precision.

\subsection{Constructing modular forms and cusp forms from the Hauptmodul}
We have seen in the previous section how the Fourier expansion of the Hauptmodul can be computed from the Belyi map. In this section we discuss how complete bases of $S_k$ and $M_k$ can be constructed from this result.

By Theorem~\ref{theorem:holomorphic_modform_polynomial} every modular function on $G$
that is holomorphic outside infinity can be written as a polynomial in the Hauptmodul
$j_G$. Since $j_G$ is a modular function (i.e., weight zero form), its derivative
$j_G'(\tau) := \frac{1}{2\pi \iu} \frac{\partial}{\partial \tau} j_G(\tau)$ is a
(weakly holomorphic) modular form of weight two. Higher weight forms can be
constructed
by computing powers of $j_G'(\tau)$ and the monomial $(j_G'(\tau))^{k/2}$ is
therefore of weight $k$. If $f$ is a holomorphic modular form of weight~$k$, then
$f(\tau)/(j_G'(\tau))^{k/2}$ is a (meromorphic) modular function which has poles at
the zeros of $j_G'(\tau)^{k/2}$, which are located at the elliptic points and cusps
other than infinity. To make this modular function holomorphic outside infinity we
cancel its poles by multiplying it with the polynomial
\begin{equation}
	B(j_G(\tau)) = B_e(j_G(\tau))\cdot B_c(j_G(\tau))\,,
\end{equation}
that is designed to cancel out all the poles up to the correct order. Because $j_G$ is a modular function on $G$, multiplying a modular form by a polynomial in $j_G$ does not destroy the modularity.

Note that $j_G'(\tau)$ has zeros of order one at the cusps that are not infinity.
Therefore, we may take
\begin{equation}
	B_c(j_G(\tau)) = \prod_{c\neq i\infty}(j_G(\tau)-j_G(c))^{\alpha_c}\,,
\end{equation}
with
\begin{equation}
	\alpha_c = k/2\,.
\end{equation}

At the elliptic points, $j_G'(\tau)$ has zeros of order $n_{e_i}-1$, where
$n_{e_i}$ denotes the order of the elliptic point which is either 2 or 3.
Following from this, we construct
\begin{equation}
	B_e(j_G(\tau)) = \prod_{e}(j_G(\tau)-j_G(e))^{\beta_e}\,,
\end{equation}
with
\begin{equation}
	\beta_e = \Big\lfloor \frac{k(n_e-1)}{2n_{e}} \Big\rfloor \,.
\end{equation}
(Note that we need to divide by the order of the elliptic point since
$(j_G(\tau)-j_G(e))$ has a zero of order $n_{e}$, see for example
\cite[pp.~227--228]{cox}.) By construction, $f(\tau)/(j'_G(\tau))^{k/2}\cdot
B(j_G(\tau))$ is a modular function that is holomorphic outside infinity and hence by
Theorem~\ref{theorem:holomorphic_modform_polynomial}
\begin{equation}
	f(\tau)/(j'_G(\tau))^{k/2}\cdot B(j_G(\tau)) = P(j_G(\tau))\,.
	\label{eq:P_j_G}
\end{equation}

We now use Eq.~\eqref{eq:P_j_G} to construct modular forms with prescribed
valuations at the cusps which can be used to construct bases of $S_k$ and $M_k$.
These constructed forms have valuations at the cusps that are equivalent to those of
a reduced row echelon basis and are therefore linearly independent.
From Eq.~\eqref{eq:P_j_G} we get that
\begin{equation}
	f(\tau) = (j'_G(\tau))^{k/2}\cdot \frac{P(j_G(\tau))}{B(j_G(\tau))}\,.
	\label{eq:f_from_polynomial}
\end{equation}
In order to get a basis of forms of $M_k$, the $i$-th form $f_i$ should have
valuation $i$ at infinity, where $i=0,1,\dots$. Note that $j_G(\tau)$ and
$j'_G(\tau)$ both have a pole of order 1 at infinity (i.e., valuation $-1$ in terms
of $q_N$). We therefore get the desired behavior at infinity by choosing
$P_i(j_G(\tau))$ to be a monomial
\begin{equation}
	P_i(j_G(\tau)) = j_G(\tau)^{\textrm{deg}(B) - k/2 - i}\,.
\end{equation}
The construction of cusp forms $f_i \in S_k$ works similarly. In this case $f_i$
should have valuation~$1$ at all cusps outside infinity and valuation $i+1$ at
infinity. We hence get
\begin{equation}
	\textrm{deg}(P_i) = \textrm{deg}(B) - k/2 - i - 1\,.
\end{equation}
In order to impose vanishing at the cusps outside infinity we simply need to multiply
by the factors $\left(j_G(\tau)-j_G(c)\right)$. Let $n(c)$ denote the number of cusps
of $G$. Then we get
\begin{equation}
	P_i(j_G(\tau)) = \prod_{c\neq i \infty}\left(j_G(\tau)-j_G(c)\right)\cdot
	j_G(\tau)^{\textrm{deg}(B) - k/2 - i - 1-(n(c)-1)}\,.
\end{equation}
\begin{example}[Constructing cusp form from Hauptmodul]
	Continuing with the example from this section, suppose that we would like to
	construct $f_0 \in S_4(G)$. By applying the result of
	Section~\ref{sec:hauptmodul_from_belyi}, we compute the $q$-expansion of the
	Hauptmodul
	\begin{equation}
		j_{G}(\tau) = q_3^{-1} + 148932u^2q_3 + 35666932u^3q_3^2 + 7392301056u^4q_3^3
		+ \dots\,.
	\end{equation}
	The space $S_4(G)$ is one-dimensional and we get
	\begin{align}
		f(\tau) &= (j'_G(\tau))^2\cdot \frac{P(j_G(\tau))}{B(j_G(\tau))}\\
		&= (j'_G(\tau))^2\cdot \frac{(j_G(\tau)+462u)}{(j_G(\tau)+462u)^2(j_G(\tau)-996u)(j_G(\tau)+516u)}\,,
	\end{align}
	which yields in the expansion
	\begin{equation}
		f(\tau) = q_3 + 18uq_3^2 - 8640u^2q_3^3 - 1823860u^3q_3^4 + \dots\,.
	\end{equation}
\end{example}
The approach presented in this section can be used to explicitly compute modular forms and cusp forms over $L$ which means that the results are rigorous. An additional advantage from a \emph{performance perspective} is that, once the Fourier expansion of the Hauptmodul has been computed, the remaining forms can be obtained without additional expensive solving or series reversion. We note however that the division of power series can be ill-conditioned when working over $\CC$ for the problems involved. For this reason it is useful to make use of \textsc{Arb}'s interval arithmetic in order to assert that the coefficients have been computed to sufficient accuracy. Once a basis of forms has been constructed, linear algebra can be used to transform the basis into reduced row echelon form.
\begin{remark}
It would be interesting to examine the practicality and effectiveness of an approach
where higher genus Newton methods (see for example
\cite{belyi_db,selander-stroembergsson}) are used to compute the curve from which the
modular forms can then be constructed.
\end{remark}

\section{Conclusion}
We have shown how modular forms on general noncongruence subgroups of moderately
large index can be computed efficiently. We are currently working on applying the
presented algorithms to create a database of modular forms and cusp forms on
noncongruence subgroups and plan to release this to the LMFDB \cite{lmfdb}. We remark
that the improved solving techniques that were presented in
Section~\ref{sec:mixed_prec_solving} should also be beneficial in the computation of
Maass cusp forms, Taylor expansions of modular forms and other examples of modular
forms at arbitrary precision arithmetic. We also hope that our demonstration of the
effectiveness of the usage of mixed-precision arithmetic in the context of arbitrary
precision arithmetic might inspire future work.

\section*{Acknowledgements}
The authors would like to thank Fredrik Strömberg for assistance in the installation of \textsc{Psage} and John Voight for useful comments.

\bibliographystyle{abbrv}
\bibliography{ms.bib}
\end{document}